\documentclass[a4paper]{article}
\title{Cayley-Menger extension of metrics\\ on affine spaces.}
\author{Ana Casimiro$^1$ \and César Rodrigo$^2$}
\date{\small%
	$^1$Department of Mathematics, NOVA SST, Universidade NOVA de Lisboa, Campus de Caparica 2829-516, Caparica, Portugal\\%
	Center for Mathematics and Applications (NOVA Math)\\
	\texttt{amc.fct.unl.pt}\\ [2ex]
	$^2$Instituto Politécnico de Setúbal, ESTSetúbal, Campus do IPS, Estefanilha 2914-761 Setúbal, Portugal\\
	Centro de Matemática, Aplicações Fundamentais e Investigação Operacional\\
	\texttt{cesar.fernandez@estsetubal.ips.pt}\\
	[2ex]%
	\today
}

\usepackage{lineno,hyperref}
\modulolinenumbers[5]
\usepackage{amsmath,amssymb,amsthm}
\usepackage[all,cmtip]{xy}
\usepackage{bbold}
\usepackage{tikz-cd}
\usepackage{tabularx}
\usepackage{makecell}

\newtheorem{define}{Definition}
\numberwithin{define}{section}

\numberwithin{example}{section}
\newtheorem{theorem}{Theorem}
\numberwithin{theorem}{section}
\newtheorem{prop}[theorem]{Proposition}
\newtheorem{corollary}[theorem]{Corollary}
\newtheorem{lemma}[theorem]{Lemma}
\newtheorem{remark}[theorem]{Remark}

\newcommand{\affi}[1]{\mathrm{#1}}
\newcommand{\bary}{\mathop{\mathrm{bary}}\nolimits}
\newcommand{\inv}{\mathop{\mathrm{inv}}\nolimits}
\newcommand{\Img}{\mathop{\mathrm{Im}}\nolimits}
\newcommand{\vect}[1]{\overrightarrow{#1}}

\newcommand{\Aff}[1]{\mathop{\mathrm{Aff}}\nolimits#1}
\newcommand{\Met}[1]{\mathop{\mathrm{Met}}#1}
\newcommand{\Hess}[1]{\mathop{\mathrm{Hess}}#1}

\newcommand{\Lin}[1]{\mathop{\mathrm{Lin}}(#1)}

\newcommand{\Sym}{\mathop{\mathrm{Sym}}\nolimits}
\newcommand{\Quad}[1]{\mathop{\mathrm{Quad}}\nolimits#1}

\newcommand{\Map}[2]{\mathop{\mathrm{Map}}(#1,#2)}
\newcommand{\diag}{\mathop{\mathrm{diag}}}
\newcommand{\CM}{\mathop{\mathrm{CM}}\nolimits}

\newcommand{\RR}{\mathbb{R}}

\newcommand{\WW}{\mathbb{W}}
\newcommand{\HH}{\mathbb{H}}
\newcommand{\one}{\mathbb{1}}
\newcommand{\zero}{\mathbb{0}}

\newcommand{\trz}{\mathop{\mathrm{tr}}\nolimits}
\newcommand{\dd}{\mathrm{d}}
\newcommand{\Id}{\mathrm{Id}}

\begin{document}
\maketitle 

\begin{abstract}
Associated to any affine space $\affi{A}$ endowed with a metric structure of arbitrary signature we consider the space of affine functionals operating on the space of quadratic functions of $A$. On this functional space we characterize a symmetric bilinear form derived from the metric structure in a functorial way. We explore the geometrical relations of the relevant objects in this new metric space. Their properties encode all characteristics known in the literature for euclidean squared distance matrices, Cayley-Menger matrices and determinants, squared distance coordinate systems, and Lie and Möbius sphere geometries. Birthing this form as Cayley-Menger product, it represents a geometrical foundation unifying results in all these areas, extending them to metric affine spaces or bundles.\newline 
\end{abstract}

\ \footnote{\textbf{Funding: }This work is funded by national funds through the FCT - Fundação Para a Ciência e a Tecnologia, I.P., under the scope of projects UIDB/00297/2020, UIDP/00297/2020 (Center for Mathematics and Applications) and UIDB/04561/2020, UIDP/04561/2020 (CMAF-CIO)}
\section{Introduction}
Affine metric geometry is a classical area in the basic training of a large variety of scientific disciplines. As a consequence few surprising facts are expected in the exploration of its fundamental objects. Squared distances between pairs of points on euclidean spaces (and also its generalization, the power associated to pairs of hyperspheres or hyperplanes) show a particular predilection to satisfy certain algebraic relations, an observation registered by Cayley \cite{cayley1841theorem}, as the necessity for the vanishing of a certain determinant, which is later identified as sufficient condition for a set of points to be immersed in a $n$-dimensional space, in the works of Menger \cite{menger1930untersuchungen}. Lachlan \cite{lachlan1886xv} treated the same order of ideas substituting squared distances between points by the so-called power corresponding to pairs of spheres (introduced by Darboux \cite{darboux1872relations}), and derived new algebraic relations for such values.  It was Cox \cite{cox1882systems} who realized that this power notion is a bilinear product when an appropriate linear structure is given on the space of circles of the plane. Power coordinates have also been used in the study of circles and spheres \cite{coolidge1916treatise}. 

In a somehow unrelated way, additional geometrical structures had been introduced by Möbius and Lie \cite{mobius1855theorie,Lie1872complexe} for the study of spheres, which can be seen as points on a certain projective space, enlarging the euclidean space with additional dimensions and by means of stereographic projections. On these projective spaces certain quadrics play a relevant role.
  
Due to its relevance for technological applications, the branch of distance geometry has attracted renewed interest in the present century. A large variety of applications are being developed and derived from the consideration of squared distance matrices and Cayley-Menger matrices, for the euclidean case \cite{alfakih2018euclidean,cao2006sensor,crippen1988distance,dokmanic2015euclidean,thomas2005revisiting,thomas2020clifford,liberti2014euclidean}. Many of these concepts and results are introduced and proven using particular, non-functorial objects, and formulated without an appropriate geometrical framework. It is not clear to which extent they are artifacts that can be created just for euclidean geometries, and how these artifacts covariate when one considers isometries between such spaces.

The somehow surprising fact, which seems to have gone unnoticed or at least never formalized, and that is presented in this article, is that all these objects, are directly related to a remarkable intrinsic bilinear form, constructed in a functorial way from any affine metric space (with arbitrary signature, and even for degenerate metrics), a bilinear form that we call Cayley-Menger bilinear form. Its functorial nature allows the exploration of its behaviour with respect to morphisms (or automorphisms) in the category of affine metric spaces. It also represents the needed concept which allows for the extension of several known results from linear geometry to the case of affine metric bundles on manifolds, bringing a new tool for the study of  (pseudo)-Riemannian manifolds.

Our paper begins in section \ref{sec2} with a presentation of affine geometry solely in terms of a barycentric operator (definition \ref{defineespacioafin}).
The affine space $\affi{A}$  and also its space of director vectors $\vect{A}$ are both immersed in a functorial way into a linear hull $\widehat{\affi{A}}$ (definition \ref{definelinearhull}). Affine functions are represented as linear functions on the linear hull, and their linear principal components are the restriction to the subspace of director vectors in the linear hull. 

Quadratic mappings and functions on affine spaces are then introduced in section \ref{sec3} in terms of its behaviour with respect to the barycentric operator (definition \ref{definequadratico}). Each quadratic function has an associated gradient covector field (definition \ref{gradientcovectorfield}) and a hessian principal component (definition \ref{definehessiancomponent}). An affine Poincare's lemma (proposition \ref{poincare}) relates affine covector fields and quadratic functions on the affine space using the gradient operator. Theorem \ref{largepropo} identifies a commutative diagram of exact sequences (\ref{hessianexactsequence}) relating symmetric bilinear forms on the linear hull with quadratic functions on the affine space.  
The hessian principal component of a quadratic function arises then as the restriction of the corresponding linear hull bilinear form to the subspace of director vectors.

The filtration of subspaces: constant functions, affine functions, quadratic functions leads to the consideration of quotient spaces: affine functions modulo constant functions (definition \ref{definegradienteafin}) is the linear space of director covectors (or constant covector fields on the affine space); quadratic functions modulo constants is the linear space of closed affine covector fields (proposition \ref{poincare}); finally quadratic functions modulo affine functions is the linear space of metrics (definition \ref{definemetricas}), which we study in section \ref{sec4}. Each quadratic function has a metric-equivalent reduced form, using homogeneization centered at a given point, or using trivialization at a given referential (proposition \ref{formulametrica}).

The determination of the homogeneous representative of the metric centered at different points $P\in\affi{A}$ gives a quadratic immersion (definition \ref{coroCMimmersiondef}) of the affine space $\affi{A}$  in the affine space of quadratic functions associated to a metric $m$. The determination of a metric by means of a quadratic function on the vector hull that vanishes at a given referential leads to the matrix representation (\ref{CMmatrix}) in proposition \ref{nondegeneracy}, which lets an insight of the relation of Cayley-Menger matrices and determinants with the interpretation of metrics as restrictions of bilinear forms on the vector hull (remark \ref{remarkEDM}).

After completion of our description of quadratic functions and metrics in terms of bilinear forms on the linear hull, section \ref{sec5} is devoted to the main result on this work, with the intrinsic introduction of a Cayley-Menger bilinear form on the quadratic hull associated to any given metric. Firstly, a quadratic hull $\widehat{\affi{A}}^2$ is introduced, together with a natural quadratic immersion of the affine space $\affi{A}$ into this space (definition \ref{define_vP}). Using the immersion any quadratic function on $\affi{A}$ can be represented as a linear function on the quadratic hull. Along this section we explore the relation of this natural immersion to Cayley-Menger mapping given in definition \ref{coroCMimmersiondef}.

Fixing a non-null metric $m\in\Met\affi{A}$ determines a quotient space of the quadratic hull, namely the space $\Aff\Quad_m\affi{A}$ of affine functionals on the space $\Quad_m\affi{A}$ of quadratic functions associated to $m$. Its dual space is precisely the linear hull $\widehat{\Quad_m\affi{A}}$ of the affine space $\Quad_m\affi{A}$, which can be seen as the linear subspace of quadratic functions whose hessian principal component is a multiple of $m$ (see lemma \ref{lemadeltav} and the following remark).

Proposition \ref{basevR} shows that any affine referential determines a natural basis on $\Aff\Quad_m\affi{A}$ and identifies its dual one. For the case of a non-degenerate metric proposition \ref{propobasedelta} determines another basis (non-dual to the former) on $\widehat{\Quad_m\affi{A}}$. This latter one is an overdetermined system of quadratic coordinates on $\affi{A}$ known as Cayley-Menger coordinates (definition \ref{CMcoordinates}). The change of basis matrix is precisely the Cayley-Menger matrix. However there exists a more fundamental aspect to this matrix. Specifically, the bilinear form $\CM_m\colon \Aff\Quad_m\affi{A}\rightarrow \widehat{\Quad_m\affi{A}}$ determined by these two basis turns out to be independent of the referential choice, and can be described in intrinsic terms with a functorial characterization. This main result of our work is stated in theorem \ref{teoremaprincipal}.

There exists, associated to each metric $m\in\Met\affi{A}$, a corresponding bilinear form $\CM_m$ (Cayley-Menger product associated to $m$) whose most elementary properties are described in proposition \ref{propCMref}, in particular it is related to Cayley-Menger matrices, to Cayley-Menger overdetermined quadratic coordinate systems, and enables the treatment and computations with $m$-quadratic functions just by using linear expressions: as shown by lemma \ref{immersproj} even though $\affi{A}$ is immersed by a quadratic mapping into $\widehat{\Quad_m\affi{A}}$ and into its dual space, the linear projector taking quotient by $\RR$ transforms these quadratic mappings into linear ones, which grants a simple (linear) mechanism to determine affine coordinates of points from the overdetermined system of quadratic coordinates.

 For example, baricentric coordinates of a point are recovered as affine combinations of the values of its Cayley-Menger coordinate (quadratic!) functions. Its inertia index is directly related to the inertia index of the metric (corollary \ref{signaturaCM}), and $\CM_m$ is non-degenerate for any non-degenerate metric $m$ (corollary \ref{nondegenerateCM}). The definition of this bilinear form on the space of $m$-quadratic functions is also intimately related (proposition \ref{RRAdual}) to theories that model hyperplanes and hyperspheres as points of a projective space, following ideas of Lie, Möbius or Pedoe \cite{Lie1872complexe,mobius1855theorie,pedoe1937representation}. However those models are focused on the euclidean case, with ad-hoc elements, like a choice of origin or stereographic projections. We may explore this question using the intrinsic Cayley-Menger product on the quadratic hull, with results to be presented in an companion paper.

In theorem \ref{CMquadricAff} and corollary \ref{invCMquadric} we also state the relation of this bilinear form with the immersion of points as evaluation operators on the space of quadratic functions (presented in proposition \ref{v_P_isquadratic}), and with the Cayley-Menger immersion (definition \ref{coroCMimmersiondef}).

In the non-degenerate case the inverse Cayley-Menger bilinear product allows a definition of (inverse) Cayley-Menger bilinear product on the linear hull of the affine space of $m$-quadratic functions. Corollary \ref{invCMquadric} shows that elements on the Cayley-Menger quadric are precisely points on the linear hull $\widehat{\Quad_m\affi{A}}$ that are isotropic with respect to this inverse Cayley-Menger bilinear product. Finally we present a theorem \ref{functorialCM} that illustrates the behavior of affine mappings with respect to the functor that associates any metric affine space $(\affi{A},m)$ to its corresponding vector space $(\Aff\Quad_m\affi{A},\CM_m)$ with Cayley-Menger bilinear form.

Main results of section \ref{sec5} are summarized in figure \ref{figura}, explained in remark \ref{diagrama}, at the end of this paper, and its computational applications are briefly illustrated in section \ref{secejemplo}. The reader is advised to check there for a specific presentation of most relevant definitions, and for the visualization of the several relations among them.

\section{Affine spaces and linear hull}\label{sec2}
We assume basic knowledge of linear and affine geometry, and of vector space duality. In the following the (dual) space of linear forms (covectors) associated to a vector space $\vect{V}$ is represented as $\vect{V}^*$, duality coupling of vectors and covectors will be represented as $\langle\,,\,\rangle$,  $f^*\colon \vect{V}^*\rightarrow \vect{E}^*$ stands for the dual morphism associated to a morphism $f\colon \vect{E}\rightarrow \vect{V}$, and the (incident) subspace of covectors vanishing on a subspace $\vect{F}\subseteq\vect{V}$ shall be represented as $\vect{F}^\circ\subseteq \vect{V}^*$. For simplicity all spaces under consideration are assumed real and finite-dimensional ones.

For convenience our approach to affine geometry will mainly rely on barycentric description of affine spaces, rather than the common approach that fixes the underlying director vector space to define the affine structure. In order to determine a (real) affine structure on some set $\affi{A}$ the only necessary element is a notion of product (or dilation), which we call barycentric operator.  Barycentric operators satisfy a certain set of axioms described in Theorem 1.1 in \cite{bertram2004linear}, or in Theorem 2.1 in \cite{buliga2010infinitesimal}. We refer to these papers for further details and assume the consequences of these axioms and all general notions and main results in this area are known to the reader:
\begin{define}\label{defineespacioafin}
We call real affine space any set $\affi{A}$ with a barycentric operator:
$$\bary\colon (\alpha,P,Q)\in\RR\times \affi{A}\times\affi{A}\mapsto \bary_\alpha(P,Q)\in\affi{A}$$
for which the following axioms hold:
		\begin{itemize}
		\item[Af1)] Multiplicative structure preservation:
		$$\bary_1(P,Q)=Q,\quad \bary_\alpha(P,P)=P,$$ $$ \bary_\alpha(P,\bary_\beta(P,Q))=\bary_{\alpha\beta}(P,Q)$$
		\item[Af2)] Each $\bary_{\alpha}(P,\cdot)$ is a morphism for each $\bary_\beta$:
		$$\bary_\alpha(P,\bary_\beta(Q,R))=\bary_\beta(\bary_\alpha(P,Q),\bary_\alpha(P,R))$$
		\item[Af3)] Barycentric condition:
		$$\bary_\alpha(P,Q)=\bary_{1-\alpha}(Q,P)$$ 
		\item[Af4)] Commutativity of translation mappings  $\bary_\alpha(P,\cdot)\circ \bary_{1/\alpha}(Q,\cdot)$:
		$$\begin{aligned}\bary_\alpha(P,\cdot)&\circ \bary_{1/\alpha}(Q,\cdot)\circ 
			\bary_\beta(R,\cdot)\circ \bary_{1/\beta}(S,\cdot)=\\&
			\bary_\beta(R,\cdot)\circ \bary_{1/\beta}(S,\cdot)\circ \bary_\alpha(P,\cdot)\circ \bary_{1/\alpha}(Q,\cdot) 
		\end{aligned}$$
	\end{itemize}
\end{define} Elements $\bary_\alpha(P,Q)$ are usually represented with classical notation as $P+\alpha\cdot \vect{PQ}$ or $(1-\alpha)\cdot P+\alpha\cdot Q$ (also called a weighed mean). 
\begin{define}
	A mapping $f\colon \affi{A}\rightarrow \affi{B}$ between affine spaces is affine if $$\bary_\alpha(f(P),f(Q))=f(\bary_\alpha(P,Q)),\quad \forall(\alpha,P,Q)\in\RR\times\affi{A}\times\affi{A}.$$
	
	The set of affine mappings between affine spaces $\affi{A}$ and $\affi{B}$ is an affine subspace of $\Map{\affi{A}}{\affi{B}}$ and will be denoted by $\Aff(\affi{A},\affi{B})$.
\end{define}

We shall call $n$-dimensional weight any column matrix  $w=[w_0\ldots w_n]^t\in M_{(n+1)1}(\RR)$ such that $\sum w_i=1$. The set of $n$-dimensional weights is an affine subspace $\WW_n\subset M_{(n+1)1}(\RR)$, and will be the basic model of ($n$-dimensional) affine space.

Any ordered sequence of $k+1$ points $\mathrm{P}=[P_0\ldots P_k]\in M_{1(k+1)}(\affi{A})$ determines an affine mapping (affine combinations of points) $w\in\WW_k\mapsto w_0P_0+w_1P_1+\ldots+w_kP_k\in\affi{A}$.

\begin{define}
	We call affine span of a set of points $P_0,\ldots, P_k\in\affi{A}$ and represent by $\langle P_0,\ldots, P_k\rangle$ the minimal affine subspace of $\affi{A}$ that contains these points. Its points are all affine combinations formed by them and weights $w=[w_0\ldots w_k]^t\in\WW_k$.
	
	When the affine space is spanned by a finite number of points, we say it is finite-dimensional. Any $\mathcal{R}=(R_0,R_1,\ldots, R_n)$ minimal ordered sequence of points  whose affine span is  $\affi{A}$ is called an affine referential of this affine space.
\end{define}
In the following we will always assume that our affine spaces are finite-dimensional. 

Affine mappings$f\colon \affi{A}\rightarrow \affi{B}$ are characterized by $f(\alpha\cdot P+\beta\cdot Q)=\alpha\cdot f(P)+\beta\cdot f(Q)$ when $P,Q\in\affi{A}$ and $\alpha+\beta=1$. In this case the property extends to arbitrary affine combinations:
$$f(w_0P_0+\ldots+w_kP_k)=w_0f(P_0)+\ldots+w_kf(P_k)$$ for arbitrary $k$-dimensional weight $w\in\WW_k$ and points $P_0,\ldots, P_k\in\affi{A}$.
In matrix form, for $P=w_0P_0+\ldots+w_kP_k$ we may write: $$f(P)=\left[Q_0\,Q_1\,\ldots \,Q_k\right]\cdot \left[\begin{matrix}w_0\\ w_1\\\vdots\\ w_k\end{matrix}\right], \qquad Q_i=f(P_i)\in\affi{B},\quad w\in\WW_k\subset M_{(k+1)1}(\RR)$$
with the obvious definition of product (in a barycentric sense), for a row vector with $k+1$ entries in $\affi{B}$ and a $k$-dimensional weight.

Considering the specific case $\affi{B}=\RR$ leads to the space of affine functions $\Aff{\affi{A}}\subset\Map{\affi{A}}{\RR}$, a vector subspace of the space $\Map{\affi{A}}{\RR}$ of real-valued functions on $\affi{A}$. Constant functions are a particular instance of affine function. In particular the constant unit function $u\in\Aff{\affi{A}}$ is affine. We have a natural linear immersion $\alpha\in\RR\mapsto \alpha u\in\Aff{\affi{A}}$.

\begin{define}\label{definegradienteafin} The quotient space $\Aff{\affi{A}}/\RR$ shall be called covector space associated to $\affi{A}$, and represented as $\vect{A}^*$.
	
	Any affine function $f$ determines $\dd f=f+\RR\in\Aff\affi{A}/\RR=\vect{A}^*$, its gradient covector, also called (linear) principal component  associated to $f$. 
\end{define}
The inclusion $\alpha\in\RR\mapsto \alpha u\in\Aff{\affi{A}}$ and the projector $f\in\Aff{\affi{A}}\mapsto \dd f\in\vect{A}^*$ determine an exact sequence of linear mappings:
\begin{equation} \label{exataafin} 0\rightarrow \RR\hookrightarrow \Aff{\affi{A}}\rightarrow \vect{A}^*\rightarrow 0.
\end{equation}

The following immersion of any affine space into a linear hull is due to Tisseron \cite{tisseron1983geometries}. The interested reader may consult some aspects of the linear hull (also called vector hull) in \cite{gracia2008vector}. For brevity we state next, without proof, its main characteristics:
\begin{define}\label{definelinearhull}
	The vector space $\widehat{\affi{A}}=(\Aff{\affi{A}})^*$ of linear operators on the space of affine functions is called linear hull associated to $\affi{A}$. Each point $P\in\affi{A}$ defines an element $z_P$ on $\widehat{\affi{A}}$, determined by $\langle z_P,f\rangle=f(P)$, for any affine function $f\in\Aff{\affi{A}}$. The element $z_P$ is called evaluation operator at $P$ on affine functions.
	
\end{define}
The affine immersion $\affi{A}\subset\widehat{\affi{A}}$ into a vector space, identifies $\Aff\affi{A}$ with $\mathop{\mathrm{Lin}}(\widehat{\affi{A}},\RR)$. In a categorical sense, any other affine immersion of this kind could be seen as a linear hull of the affine space. This is the case of the vector space referred to as ``universal space'' $\widehat{X}$ associated to an affine space $X$ in Berger's Geometry \cite{berger2009geometry} (sections 3.1 and 3.2). We prefer Tisseron's presentation, for its intrinsic nature.

The mapping $P\in\affi{A}\mapsto z_P\in\widehat{\affi{A}}$ is an affine immersion that determines an identification of $\affi{A}$ as an affine subspace of its linear hull:
	\begin{equation}\label{formula0}\affi{A}\simeq \{z\in\widehat{\affi{A}}\,\colon\, \langle z,u\rangle=1\}\subset \widehat{\affi{A}}.
	\end{equation}
A second statement relates barycentric referentials on any affine space with linear referentials on its linear hull:
\begin{prop}\label{base_zP}
An ordered sequence of points $\mathcal{R}=(R_0,\ldots,R_n)$ is an affine referential on  $\affi{A}$ if and only if $z_{R_0},\ldots,z_{R_n}$ is a basis of the vector space $\widehat{\affi{A}}$. Its dual basis is a sequence of functions $w_0,\ldots, w_n\in\Aff{\affi{A}}$ called barycentric coordinate functions associated to the referential.
\end{prop}
Fixing a referential $\mathcal{R}=(R_0,\ldots,R_n)$ we may represent any vector $z\in\widehat{\affi{A}}$ as a linear combination $x_0z_{R_0}+x_1z_{R_1}+\ldots +x_nz_{R_n}$. The column vector $x=[x_0\ldots x_n]^t\in M_{(n+1)1}(\RR)$ is called linear coordinate vector (or briefly linear coordinates) associated to $z\in\widehat{\affi{A}}$  with respect to the basis (linear referential) $z_{R_0},\ldots,z_{R_n}$, or also with respect to the affine referential $\mathcal{R}$.
The coordinate functions $x_0,\ldots,x_n$ are linear on $\widehat{\affi{A}}$ and restrict as affine functions $w_i$ on $\affi{A}\subset\widehat{\affi{A}}$. In fact these affine functions $w_0,\ldots,w_n$ for a basis on the vector space $\Aff{\affi{A}}$. 

Barycentric coordinates are a system of $n+1$ affine functions on $\affi{A}$, and they are overdetermined: knowledge of the value of $n$ of these functions on a point implies the knowledge of the value of the remaining coordinate function. In this paper we shall deal with different notions of systems of coordinates. For precision of language it is convenient to clearly state its meaning:
\begin{define}\label{odsystem}
	We call over-determined system of coordinates on the affine space $\affi{A}$ (briefly ``o.d.system'') any (non-necessarily affine) injective mapping $x\colon \affi{A}\rightarrow \RR^m$. Its component $x_i\in\Map{\affi{A}}{\RR}$ is called $i$-th coordinate function of the system. When $x$ is bijective we will talk simply of a system of coordinates.
	
	For a fixed o.d.system of coordinates $x\colon \affi{A}\rightarrow \RR^m$, each point $P\in\affi{A}$ is univocally determined by some element $x(P)=[x_1(P)\ldots x_m(P)]^t\in M_{m1}(\RR)$ which is called the coordinate vector associated to this point for the given o.d.system of coordinates.
\end{define}
When talking of a system of affine or quadratic coordinates, we mean that the mapping $x$ is correspondingly affine, or quadratic (see definition \ref{definequadratico}). In the case that $\affi{A}$ is endowed with a vector space structure and $x$ is a linear mapping, we say that it defines a system of linear coordinates.

\begin{remark}
	Each basis $w_0,\ldots, w_n\in\Aff{\affi{A}}$ for which $w_0+w_1+\ldots+w_n=u$ holds, has a dual basis $z_0,z_1,\ldots, z_n$ that belongs to $\affi{A}\subseteq \widehat{\affi{A}}$, indeed:
	$$\langle z_i,u\rangle=\langle z_i,w_0+\ldots+w_n\rangle=\sum_j \langle z_i,w_j\rangle=1\Rightarrow z_i\in \affi{A}\subseteq \widehat{\affi{A}}.$$
	When the injective mapping $x=(w_0,\ldots, w_n)$ takes values in the space of $n$-dimensional weights it is called an o.d.system of barycentric coordinates. Its components are affine functions that determine the barycentric coordinates of any point, with respect to some affine referential.
\end{remark}
Any o.d.system of barycentric coordinates can be seen as a system of linear coordinates on $\widehat{\affi{A}}$, or as an o.d.system of affine coordinates on $\affi{A}$ (only elements with $x_0+\ldots+x_n=1$ represent a point on $\affi{A}$).

\begin{define}\label{definec}
Linear coordinates $(c_0,\ldots, c_n)$ of an affine function $f\in\Aff{\affi{A}}$ in the o.d.system of barycentric coordinates $w_0,\ldots, w_n\in\Aff{\affi{A}}$ are called linear coefficients associated to $f$ with respect to the corresponding referential. They determine a row vector $[c_0\ldots c_n]\in M_{1(n+1)}(\RR)$, called coefficient row vector (or briefly coefficients) associated to the affine function, in the corresponding affine referential $\mathcal{R}=(R_0,\ldots,R_n)$.
\end{define}
For any function $f\in\Aff{\affi{A}}$ and point $Q\in\affi{A}$ represented by coefficients $c=[c_0\ldots c_n]\in M_{1(n+1)}(\RR)$, and an $n$-dimensional weight $w=[w_0\ldots w_n]^t$ respectively, the value $f(Q)$ is simply the matrix product $c\cdot w$.

The unit function $u\in\Aff{\affi{A}}$ is represented, in terms of its coefficients with respect to any affine referential, by the unit row vector $[1\ldots 1]=\one_n\in M_{1\times(n+1)}(\RR)$.

Using the immersion $\affi{A}\subset \widehat{\affi{A}}$, points $Q\in\affi{A}$ are represented by a coordinate vector $q\in M_{(n+1)1}(\RR)$ in a basis $v_{R_0},\ldots,v_{R_n}$, for which $\one_n\cdot q=1$, hence determining a $n$-dimensional weight  $q\in\WW_n$, which we call barycentric coordinate vector associated to $Q$ with respect to the referential $\mathcal{R}=(R_0,\ldots, R_n)$.

Elements in the covector space $\vect{A}^*=\Aff{\affi{A}}/\RR$ can be identified with equivalence classes $c+\RR\one_n$ of $(n+1)$-row vectors $c$, with respect to the addition of multiples of the unit row vector $\one_n$. That is, a choice of affine referential determines a linear isomorphism $\vect{A}^*\simeq M_{1(n+1)}(\RR)/\langle \one_n\rangle$.

\begin{define} Linear operators on the covector space $\vect{A}^*=\Aff{\affi{A}}/\langle u\rangle$ are called director vectors of the affine space $\affi{A}$. This space of linear operators (dual space of $\vect{A}^*$) is represented as $\vect{A}$ and called director vector space associated to $\affi{A}$.\end{define}
Taking into account (\ref{formula0}) and the exact sequence dual to (\ref{exataafin}), $\vect{A}$ can be seen as a linear subspace of $(\Aff{\affi{A}})^*$, the space of elements that are incident with $u\in\Aff{\affi{A}}$.
$$\vect{A}=\langle u\rangle^\circ =\{z\in\widehat{\affi{A}}\,\colon\, \langle z,u\rangle=0\}\subset \widehat{\affi{A}}=(\Aff{\affi{A}})^*.$$

\begin{define} Denote $\vect{PQ}=z_Q-z_P\in \widehat{\affi{A}}$. This element is incident with $u\in\Aff{\affi{A}}$, hence it belongs to the director vector space  $\vect{A}$ associated to $\affi{A}$, and is called the director vector associated to the pair of points $P,Q\in\affi{A}$.
\end{define}

When $\vect{x}\in\vect{A}\subset\widehat{\affi{A}}$ has coordinate vector $x=[x_0\ldots x_n]^t$ in the referential $\mathcal{R}=(R_0,\ldots, R_n)$, there holds $x_0+\ldots+x_n=0$. The sequence $(x_0,\ldots, x_n)$ belongs to the subspace $\langle \one_n\rangle^\circ\subseteq M_{(n+1)1}(\RR)$. We call this the space of hollow weights. \begin{equation}\label{hollowweights} \HH_n=\langle \one_n\rangle^\circ=\left\{[x_0\ldots x_n]^t\in M_{(n+1)1}(\RR)\,\colon \sum x_i=0\right\}.\end{equation}
Removing the first coordinate of the coordinate vector associated to $\vect{x}$ we get $(x_1,\ldots, x_n)$, the linear coordinates  of $\vect{x}$ with respect to the basis $\vect{R_0R_1},$ $\ldots,$ $\vect{R_0R_n}$ on $\vect{A}$.

Recall that the composition of affine mappings is again affine. Any affine mapping $f\in\Aff{(\affi{A},\affi{B})}$ determines, by composition, a linear mapping $\Aff{\affi{B}}\rightarrow \Aff{\affi{A}}$, and the corresponding dual linear mapping $\widehat{f}\colon \widehat{\affi{A}}\rightarrow \widehat{\affi{B}}$ characterized by the condition $\widehat{f}(z_P)=z_{f(P)}$, for each element $P$ on the affine subspace $\affi{A}\subset \widehat{\affi{A}}$. Considering the behavior $u_B\circ f=u_A$ for the unit functions, we conclude that $\widehat{f}$ also takes elements of the vector subspace $\vect{A}\subset \widehat{\affi{A}}$ into elements of $\vect{B}\subset \widehat{\affi{B}}$. The restriction of $\widehat{f}$ to $\vect{A}$ is represented as $\dd f\colon \vect{A}\rightarrow \vect{B}$ and called (linear) principal component associated to the affine mapping $f$. For the particular case $\affi{B}=\RR$, this notion coincides with the previously defined gradient covector $\dd f\in\vect{A}^*$ associated to functions on $\affi{A}$. Observe that:
$$\dd f(\overrightarrow{PQ})=\widehat{f}(\overrightarrow{PQ})=\widehat{f}(z_Q-z_P)=\widehat{f}(z_Q)-\widehat{f}(z_P)=z_{f(Q)}-z_{f(P)}=\overrightarrow{f(P)f(Q)}.$$

Affine mappings $f$ between two affine spaces $\affi{A},\affi{B}$ are identified with linear mappings $\widehat{f}$ between the corresponding linear hulls $\widehat{\affi{A}},\widehat{\affi{B}}$, such that $\widehat{f}(\affi{A})\subseteq \affi{B}\subset \widehat{\affi{B}}$. Using the o.d.systems of barycentric coordinates with respect to a referential $\mathcal{R}=(R_0,\ldots, R_n)$ on $\affi{A}$ and with respect to another referential $\mathcal{S}=(S_0,\ldots, S_m)$ on $\affi{B}$, the affine space of affine mappings $f\colon \affi{A}\rightarrow \affi{B}$ is identified with the affine subspace $\WW_{m,n}$ of matrices $C\in M_{(m+1)(n+1)}(\RR)$ such that $\one_m\cdot C=\one_n$ (called weight matrices). If $P\in\affi{A}$ has barycentric coordinate vector $p\in\WW_n$ then $f(P)\in\affi{B}$ has barycentric coordinate vector $C\cdot p\in\WW_m$. The linear mapping $\widehat{f}$ associated to $f$ is described in the linear basis $z_{R_0},\ldots,z_{R_n}$ of $\widehat{\affi{A}}$ and $z_{S_0},\ldots, z_{S_m}$ of $\widehat{\affi{B}}$ using matrix product with $C$. The linear mapping $\vect{f}$ induced on $\vect{A}$ is simply the restriction of this matrix product to the subspace of hollow weights $\HH_n$ (\ref{hollowweights}). Constant affine mappings are represented by matrices $C=c\cdot \one_n$, for any choice of $c\in\WW_m$. Composition of affine mappings is represented by matrix product.

\begin{remark}
The inclusion of any affine subspace $\affi{F}\subseteq \affi{A}$, induces a projector $\Aff{\affi{A}}\rightarrow \Aff{\affi{F}}$ such that $u_{\affi{A}}$ projects to $u_{\affi{F}}$, and an immersion $\widehat{\affi{F}}\subseteq \widehat{\affi{A}}$, such that $\vect{F}=\langle u_{\affi{F}}\rangle^\circ$ is a subspace of $\vect{A}=\langle u_{\affi{A}}\rangle^\circ$. The director vector space of an affine subspace $\affi{F}\subseteq \affi{A}$ is a linear subspace $\vect{F}\subseteq\vect{A}$.

Affine subspaces $\affi{F}\subseteq \affi{A}$ that have a specific director subspace $\vect{F}\subseteq\vect{A}$ can be written in the form $P+\vect{F}$. This set of affine subspaces with a fixed director subspace $\vect{F}\subseteq \vect{A}$ is again an affine space with the obvious definition of barycenter:
$$(1-\alpha)(P+\vect{F})+\alpha(Q+\vect{F})=(1-\alpha)P+\alpha Q+\vect{F},$$
and is called quotient affine space $\affi{A}/\vect{F}$. Its director vector space is the quotient vector space $\vect{A}/\vect{F}$.
\end{remark}

\section{Quadratic functions on affine spaces}\label{sec3}
A characteristic aspect of quadratic mappings is that all its values on a line can be interpolated from the values at two points and its middle point. On the other hand, on any affine space the condition of being quadratic can be expressed as being quadratic on each of its lines, hence the following definition:
\begin{define}\label{definequadratico}
	A mapping between affine spaces $\delta\colon\affi{A}\rightarrow \affi{B}$ is called a quadratic mapping if for any pair of points $P,Q\in\affi{A}$ and 1-dimensional weight $[\alpha\,\beta]^t\in\WW_1$, there holds:
	\begin{equation}\label{quadraticity}\delta(\alpha P+\beta Q)=\alpha(\alpha-\beta)\cdot \delta (P)+\beta(\beta-\alpha)\cdot  \delta(Q)+4\alpha\beta\cdot \delta\left(\frac{P+Q}{2}\right).\end{equation}
\end{define}
This formula is given only when $\alpha+\beta=1$, in which case $\alpha(\alpha-\beta)+\beta(\beta-\alpha)+4\alpha\beta=\alpha^2+\beta^2+2\alpha\beta=(\alpha+\beta)^2=1$ and this affine combination of 3 points with 2-dimensional weight $[\alpha^2-\alpha\beta\, \beta^2-\alpha\beta\, 4\alpha\beta]^t\in\WW_2$ makes sense.
	
	From this definition certain properties are easy to prove:
		\begin{itemize}
			\item Affine mappings are quadratic mappings. A quadratic mapping $\delta$ is affine if and only if $\delta\left(\frac{P+Q}{2}\right)=\frac12\delta(P)+\frac12\delta(Q)$.
			\item Composition of an affine mapping with a quadratic mapping is quadratic.
			\item Any affine combination of quadratic mappings is quadratic: quadratic mappings between affine spaces are an affine subspace $\Quad({\affi{A}},{\affi{B}})\subseteq \Map{\affi{A}}{\affi{B}}$.
		\end{itemize}
In the finite-dimensional case it is convenient to give a coordinate matrix representation of quadratic mappings. In the presence of an affine referential this representation is obtained by degree 2 homogeneous polynomials in the barycentric coordinates. Our coordinate-free notion of quadratic mapping leads to the following coordinate characterization. We provide a proof in the appendix of this work:
\begin{theorem}\label{thmmatrixquad}
	Consider two affine spaces $\affi{A}$ and $\affi{B}$. Fix an affine referential $\mathcal{R}=(R_0,R_1,\ldots,R_n)$ of $\affi{A}$ and its midpoints $R_{ij}=\frac{R_i+R_j}{2}$. For any given $\binom{n+2}{2}$ points $Q_{ij}\in\affi{B}$ (where $Q_{ij}=Q_{ji}$, $0\leq i,j\leq n$), there exists a unique quadratic  mapping $\delta\colon \affi{A}\rightarrow \affi{B}$ such that $\delta(R_{ij})=Q_{ij}$. This quadratic  mapping can be given as:
	\begin{equation}\label{matrizDelta} \delta(P)=\delta(p_0 R_0+\ldots+p_nR_n)=\sum_{i,j} p_ip_j \Delta_{ij},\qquad \Delta_{ij}=2Q_{ij}-\frac12 Q_{ii}-\frac12 Q_{jj}\in\affi{B}.\end{equation}
\end{theorem}
Observe that $2-\frac12-\frac12=1$ and that $\sum p_ip_j=1$ for any $p\in\WW_n$, hence all formulas in the theorem are valid affine combinations.

In other words, if we fix a referential $\mathcal{R}$ on $\affi{A}$, $\Quad(\affi{A},\affi{B})$ can be identified with the affine space of squared $\affi{B}$-valued symmetric matrices $\Delta\in \Sym_{n+1}(\affi{B})$. If $\delta\in\Quad(\affi{A},\affi{B})$ is represented by a $\affi{B}$-valued symmetric matrix $\Delta$, the image of $P\in\affi{A}$ can be written, with the obvious $B$-valued matrix notation, as:
$$\delta(P)=p^t\cdot \Delta\cdot p$$
where $p\in\WW_n$ is the barycentric coordinate vector associated to $P$ in the referential $\mathcal{R}=(R_0,\ldots, R_n)$, and $\Delta_{ij}=2\delta\left(\frac{R_i+R_j}{2}\right)-\frac12 \delta(R_i)-\frac12 \delta(R_j)\in\affi{B}$.

We shall now focus on quadratic functions ($\RR$-valued mappings), which represent a vector subspace $\Quad\affi{A}\subseteq\Map{\affi{A}}{\RR}$  of the space of functions on $\affi{A}$.
\begin{remark}\label{representaafines}
	Affine functions are a particular case of quadratic functions. If $c\in M_{1(n+1)}(\RR)$ is the coefficient row vector associated to some affine function $f\in\Aff{\affi{A}}$ with respect to some referential $\mathcal{R}$ (see definition \ref{definec}), then theorem \ref{thmmatrixquad} gives the alternative representation $\delta(Q)=q^t\cdot \Delta\cdot q$  with the symmetric matrix $\Delta=\frac12(\one^t\cdot c+c^t\cdot \one)$ (where $q$ stands for the barycentric coordinates of $Q$ with respect to the given referential).
\end{remark}

\begin{define}\label{defconvex}
	We say a quadratic function $\delta\in\Quad\affi{A}$ is convex if for any pair of different points $P\neq Q\in\affi{A}$ and any strictly positive weight $w=[\alpha\,\beta]\in\WW_1$ (hence with $\alpha+\beta=1$, $\alpha,\beta>0$), there holds a strict inequality:
	$$\delta(\alpha\cdot P+\beta\cdot Q) < \alpha\cdot \delta(P)+\beta\cdot \delta(Q).$$
\end{define}
Convexity is preserved when we add an affine function to a given quadratic function.

\begin{remark}\label{quadraticfrombilinear} On any vector space $\vect{V}$ one may consider the space of bilinear forms $\Lin{\vect{V},\Lin{\vect{V},\RR}}=\Lin{\vect{V},\vect{V}^*}$. For each bilinear form $g$ there exists an adjoint bilinear form $g^*$ defined by $\langle \vect{x},g^*(\vect{y})\rangle=\langle \vect{y},g(\vect{x})\rangle$. Self-adjoint bilinear forms (also called symmetric) constitute a subspace $S^2\vect{V}^*$, and any of its elements can be seen as a symmetric scalar product on $\vect{V}$, which determines a quadratic function $\delta(\vect{x})=\frac12 g(\vect{x},\vect{x})=\frac12\langle \vect{x},g(\vect{x})\rangle$ that vanishes at the zero vector and is even ($\delta(\vect{x})=\delta(-\vect{x})$).

By choosing a basis $\vect{v_1},\ldots, \vect{v_n}$ on $\vect{V}$, there exists an isomorphism of the vector space $S^2\vect{V}^*$ of symmetric bilinear forms on $\vect{V}$ with the vector space of symmetric $n\times n$ matrices. Each symmetric matrix $G$ determines a symmetric bilinear form and a quadratic function  given in linear coordinates as:
	$$\langle\vect{x},g(\vect{y})\rangle=x^t\cdot G\cdot y,\qquad \delta(\vect{x})=\frac12 x^t\cdot G\cdot x$$
	where $x,y\in M_{n1}(\RR)$ are the coordinate vectors associated to the vectors $\vect{x},\vect{y}$ in the given basis, and conversely, each symmetric bilinear form $g$ determines the values  $G_{ij}=g(\vect{v_i},\vect{v_j})$, components of the so-called Gram matrix associated to the bilinear form in the given basis.
	
	A symmetric matrix $G$ is called positive definite when $x^t\cdot G\cdot x>0$ on any non-zero column vector $x\neq 0$, which is equivalent to state that the associated quadratic function $\delta$ is convex.
\end{remark}

\begin{define}\label{simetrico}
	Consider a point $P\in\affi{A}$ on an affine space $\affi{A}$ and $\inv_P\colon Q\in\affi{A}\mapsto 2P-Q\in\affi{A}$, the inversion centered at this point $P$.  Quadratic functions $\delta$ such that $\delta\circ\inv_P=\delta$ are called even with respect to the center $P$.We denote $\Quad^h_P\affi{A}\subseteq \Quad\affi{A}$, and call subspace of homogeneous quadratic functions with respect to $P$ the set of all quadratic functions that vanish at $P$ and are even with respect to this point. 
\end{define}
\begin{define}
	For any quadratic function $\delta\in\Quad{\affi{A}}$ we call homogenization of $\delta$ with respect to $P$ the element $\delta^h_P\in\Quad^h_P\affi{A}$ defined by:
	$$\delta^h_P(Q)=\frac12\left(\delta(Q)+\delta(\inv_PQ)\right)-\delta(P).$$
\end{define}
It is clear from this definition that the homogenization is a linear retraction from the vector space $\Quad\affi{A}$ to $\Quad^h_P\affi{A}$. Using the quadraticity condition (\ref{quadraticity}):
$$\delta(\inv_PQ)=\delta(2P-Q)=6\delta(P)+3\delta(Q)-8\delta\left(\frac{P+Q}{2}\right)$$hence:
$$\delta^h_P(Q)=2\delta(P)+2\delta(Q)-4\delta\left(\frac{P+Q}{2}\right)$$
\begin{remark}
	Fix an affine referential $\mathcal{R}$ on the affine space $\affi{A}$. If $\delta\colon\affi{A}\rightarrow \RR$ is determined using theorem \ref{thmmatrixquad} as $\delta(Q)=q^t\cdot \Delta\cdot q$ for some symmetric matrix $\Delta$ (where $q\in\WW_n$ is the barycentric coordinate vector associated to  $Q$ with respect to some affine referential) then for any point $P\in\affi{A}$ with barycentric coordinate vector $p\in\WW_n$ the quadratic function $\delta^h_P$ is determined in the same referential by \begin{equation}\label{qpDelta}
		\delta^h_P(Q)=2p^t\Delta p+2q^t\Delta q-4\left(\frac{p+q}{2}\right)^t\Delta \left(\frac{p+q}{2}\right)=(q-p)^t\Delta(q-p).
	\end{equation}
	This may also be written as :
	$$\delta^h_P(Q)=q^t(\Id-p\cdot\one)^t\cdot \Delta\cdot(\Id-p\cdot\one)q.$$
	
	 Hence  the quadratic function $\delta^h_P$ is represented in the referential by:
	 \begin{equation}\label{formula4b} \Delta^h_p=(\Id-p\cdot\one)^t\cdot \Delta\cdot(\Id-p\cdot\one).\end{equation}
	 As a consequence, the function $\delta-\delta_P^h$ is affine and can be written (see remark \ref{representaafines}) as $\frac12(\one^t\cdot c+c^t\cdot \one)$, for the following coefficient row vector:
	\begin{equation}\label{coefdeltap}
		\delta-\delta^h_P\in\Aff\affi{A} \text{ has coefficient row vector }c=2p^t\Delta-p^t\Delta p\one
	\end{equation}

	This determines a decomposition: $$\Quad{\affi{A}}=\Aff\affi{A}\oplus \Quad_P^h{\affi{A}}$$
	with $\delta\mapsto(\delta-\delta^h_P,\delta^h_P)$ the natural splitting morphism.
	
	An alternative decomposition is given when we consider $\Quad{\affi{A}}=\Aff\affi{A}\oplus \Quad^0_{\mathcal{R}}{\affi{A}}$ where $\Quad^0_{\mathcal{R}}\affi{A}$ stands for the space of quadratic functions that vanish at all points of the referential $\mathcal{R}$. For a quadratic function $\delta$, its component $\delta^0_{\mathcal{R}}\in \Quad^0_{\mathcal{R}}\affi{A}$ is simply the difference of $\delta$ with the unique affine function $f$ taking values $f(R_i)=\delta(R_i)$. This affine function is characterized by a row vector $c=\diag \Delta$, where the row coefficient vector $\diag\Delta$ is the diagonal of $\Delta$. Hence following remark \ref{representaafines} the quadratic function $\delta^0_{\mathcal{R}}=\Delta-f$ is represented by:
	\begin{equation}\label{coefdelta0} \Delta^0_{\mathcal{R}}=\Delta-\frac12\left( (\diag \Delta)^t\one+\one^t(\diag\Delta)\right)
	\end{equation}

\end{remark}
	
\begin{define}\label{gradientcovectorfield}
	Consider a point $P\in\affi{A}$ and a quadratic function $\delta\in\Quad\affi{A}$ on any affine space $\affi{A}$. The affine function $\delta-\delta^h_P\in\Aff{\affi{A}}$ determines an element $\nabla_P\delta\in\Aff{\affi{A}}/\RR=\vect{A}^*$, called gradient covector associated to the quadratic function $\delta$ at the point $P$. We call gradient covector field associated to $\delta$, the affine mapping:
	$$\nabla \delta\colon P\in \affi{A}\mapsto \nabla_P\delta\in\vect{A}^*.$$
\end{define}
Observe that affine functions $\delta\in\Aff{\affi{A}}\subset\Quad{\affi{A}}$ have constant gradient covector field given by definition \ref{gradientcovectorfield}, and its (constant) value coincides with the notion in definition \ref{definegradienteafin} for the gradient covector associated to affine functions.

As we see in coordinate representation (\ref{coefdeltap}), even though $P\in\affi{A}\mapsto \delta-\delta_P^h\in\Aff\affi{A}$ is a quadratic mapping (there is a quadratic dependence on the chosen $P$), its linear principal component (we mean, the induced element modulo $\one$) depends on $P$ only in an affine fashion:   $\nabla\delta\colon P\mapsto \nabla_P\delta$ is an affine mapping on $\affi{A}$, hence an affine covector field. Moreover $\nabla\colon\Quad\affi{A}\rightarrow \Aff(\affi{A},\vect{\affi{A}}^*)$ is a linear mapping. 
\begin{define}\label{definehessiancomponent}
	We call principal component associated to an affine covector field $\omega\in \Aff(\affi{A},\vect{\affi{A}}^*)$ the associated (not necessarily symmetric) bilinear form $\vect\omega\in\Lin{\vect{\affi{A}},\vect{\affi{A}}^*}$ determined by its linear principal component (seen as mapping on the affine space $\affi{A})$. For any $\vect{x},\vect{y}\in\vect{A}$:
	\begin{equation*}
		\langle \vect{x},\vect{\omega}(\vect{y})\rangle=
		\langle \omega_{P+\vect{y}},\vect{x}\rangle-\langle \omega_{P},\vect{x}\rangle\quad(\text{arbitrary choice of }P\in\affi{A}).
	\end{equation*}
 We call hessian principal component $\Hess\delta$ associated to a quadratic function $\delta\in\Quad\affi{A}$ the principal component $\vect{\nabla\delta}$ associated to its gradient covector field $\omega=\nabla\delta$. 
\end{define}

\begin{lemma}\label{lemahesshomogeneo}
	For any quadratic function $\delta\in\Quad\affi{A}$ and for its hessian principal component $\Hess\delta\in\Lin{\vect{\affi{A}},\Lin{\vect{\affi{A}},\RR}}$ there holds:
	$$\delta^h_P(Q)=\frac12\langle \Hess\delta(\vect{PQ}),\vect{PQ}\rangle.$$
\end{lemma}
\begin{proof}
	It suffices to consider the representation of $\delta$ by some symmetric matrix $\Delta$ in some referential $\mathcal{R}$. If points $P,Q$ are represented by its barycentric coordinate vectors $p,q$, if affine functions are represented by its coefficient row vectors $c\in M_{1(n+1)}(\RR)$, if vectors $\vect{x},\vect{y}$ are represented by its coordinate hollow vectors $x,y\in\HH_n$ (vectors that are incident to $\one$), and if covectors are represented by row vectors modulo $\one$, then our previous computations show that all elements defined by $\delta$ depend on the matrix $\Delta$ as:
	$$\delta^h_P(Q)=(q-p)^t\Delta(q-p),\qquad p,q\in\WW_n \text{ (from (\ref{qpDelta}))},$$
	$$\nabla\delta\colon p\mapsto 2p^t\Delta\,(\text{mod }\one),\qquad p\in\WW_n \text{ (from (\ref{coefdeltap}))},$$
	which leads to
	\begin{equation}\label{coordrephess}
	\begin{aligned}		&\Hess\delta\colon \vect{x}\mapsto 2x^t\Delta\,(\text{mod }\one),\quad & \vect{x},\vect{y}\in\vect{A} \text{ director vectors},\\
	\frac12&\langle\Hess\delta(\vect{x}),\vect{y}\rangle=x^t\Delta y,\quad& x,y\in\HH_n\text{ coordinate hollow vectors}.\end{aligned}
	\end{equation}
	
	Observing that the vector $\vect{PQ}$ is represented by the hollow coordinate vector $x=q-p$ completes our proof.
\end{proof}
\begin{define}	
	We say an affine covector field $\omega\in\Aff(\affi{A},\vect{\affi{A}}^*)$ is a symmetric (or ``closed'') affine covector field if its principal component $\vect{\omega}$ is self-adjoint:	\begin{equation}\label{exteriorderivative} \langle \omega_{P+\vect{y}},\vect{x}\rangle-\langle \omega_{P},\vect{x}\rangle=\langle \omega_{P+\vect{x}},\vect{y}\rangle-\langle \omega_{P},\vect{y}\rangle, \quad \vect{x},\vect{y}\in\vect{A}.\end{equation}
	 Symmetric (closed) affine covector fields form a linear subspace $\Aff^S(\affi{A},\vect{\affi{A}}^*)$. In particular the gradient covector fields associated to quadratic functions belong to this subspace.
\end{define}
\begin{remark}\label{derivadaexterior}
	For any affine covector field $\omega\colon\affi{A}\rightarrow \vect{A}^*$, its exterior differential at any point $P\in\affi{A}$ is an alternating bilinear form given on tangent vectors $\vec{x},\vec{y}\in T_P\affi{A}\simeq\vect{A}$ at any point $P\in\affi{A}$ by:
	$$\begin{aligned}\dd^{ext}_P\omega(\vect{x},\vect{y})&=\langle \nabla_P(\omega(\vect{x})),\vect{y} \rangle - \langle\nabla_P (\omega(\vect{y})) ,\vect{x} \rangle=\\
	&=\langle \omega_{P+\vect{y}},\vect{x}\rangle-\langle \omega_P,\vect{x}\rangle-\langle \omega_{P+\vect{x}},\vect{y}\rangle +\langle \omega_P,\vect{y}\rangle\end{aligned}$$ 
	which justifies the interpretation of (\ref{exteriorderivative}) as a notion of closedness of affine covector fields.
\end{remark}
\begin{prop}[Affine Poincare's Lemma]\label{poincare}
	The nullspace of the gradient linear mapping $\nabla\colon \Quad\affi{A}\rightarrow \Aff(\affi{A},\vect{\affi{A}}^*)$  is $\RR\subset\Quad\affi{A}$, and its image is the subspace $\Aff^S(\affi{A},\vect{\affi{A}}^*)$ of affine covector fields $\omega\colon\Aff(\affi{A},\vect{\affi{A}}^*)$ with self-adjoint principal component.
\end{prop}
\begin{proof}
	Clearly, by definition \ref{gradientcovectorfield}, the gradient of a constant function $\delta\in\RR\subset\Quad\affi{A}$ vanishes. On the other hand  a quadratic function such that $\nabla\delta=0$ has $\delta-\delta^h_R$ a constant function for each point $R\in\affi{A}$. The function $\delta$ is then even with respect to any point $R$ and taking into account that $P$ is the symmetric point of $Q$ with respect to $R=\frac12P+\frac12Q$, we conclude $\delta(P)=\delta(Q)$ for any pair of points $P,Q\in\affi{A}$. This completes our characterization of the nullspace.
	
	Regarding the image, $\Hess\delta$ has been proven to be symmetric in (\ref{coordrephess}). 
	We have then an exact sequence of linear mappings:
	$$0\rightarrow\RR\hookrightarrow \Quad\affi{A}\stackrel{\nabla}{\longrightarrow} \Img\nabla\rightarrow 0,$$
	where $\Img\nabla\subseteq\Aff^S(\affi{A},\vect{\affi{A}}^*)$ is a linear subspace.
	
	From $\dim\affi{A}=n$ we know $\dim\Quad\affi{A}=\frac{(n+2)(n+1)}{2}$, hence $\dim \Img\nabla=\frac{n^2+3n+2}{2}-1$.
	
	Observe that $\Aff^S(\affi{A},\vect{\affi{A}}^*)$ is the inverse image of $S^2\vect{\affi{A}}^*\subset \Lin{\vect{\affi{A}},\vect{\affi{A}}^*}$
	by the mapping $\dd\colon \omega\in \Aff{(\affi{A},\vect{\affi{A}}^*)}\mapsto \vect{\omega}\in\Lin{\vect{\affi{A}},\vect{\affi{A}}^*}$. This is a surjective linear mapping whose nullspace is the space of constant $\vect{\affi{A}}^*$-valued fields. Therefore:
	$$\dim\Aff^S(\affi{A},\vect{\affi{A}}^*)=\dim\vect{\affi{A}}^*+\dim S^2\vect{\affi{A}}^*=n+\frac{n(n+1)}{2}.$$
	By dimension we must conclude that $\Img\nabla=\Aff^S(\affi{A},\vect{\affi{A}}^*)$ thus completing our proof.
\end{proof}
Following remark \ref{derivadaexterior}  the previous proposition is a formulation in affine geometry of Poincaré's Lemma. Namely, condition (\ref{exteriorderivative}) could be stated as the vanishing of the exterior derivative $(\dd^{ext}\omega)_P(\vect{x},\vect{y})$ if we take $\omega$ as a smooth 1-form, $\dd^{ext}$ the exterior differential of differential forms, and $\vect{x},\vect{y}$ as elements of the tangent space $T_P\affi{A}$. We are stating that any closed affine covector field is the differential of some quadratic function, and that the only functions whose differential vanishes are constant functions.

As we see next, in the same manner as (non-homogeneous) affine functions $f$ on $\affi{A}$ are determined by linear forms $\widehat{f}\in \widehat{\affi{A}}^*$ on its linear hull $\widehat{\affi{A}}$, also (non-homogeneous) quadratic functions $\delta$ on $\affi{A}$ are determined by symmetric bilinear forms $\widehat{g}\in S^2\widehat{\affi{A}}^*=S^2\Aff{\affi{A}}$ on the linear hull.

\begin{theorem}\label{largepropo}
Consider the linear hull $\widehat{\affi{A}}=(\Aff\affi{A})^*$ of some affine space $\affi{A}$. The mapping that determines the quadratic function $\delta(P)=\frac12\widehat{g}(z_P,z_P)$ on $\affi{A}$ associated to any symmetric bilinear form $\widehat{g}\in S^2\widehat{\affi{A}}^*=S^2\Aff\affi{A}$ is a linear isomorphism closing a commutative diagram of linear morphisms, where each of the rows is an exact sequence:
	\begin{equation}\label{hessianexactsequence}
		\begin{tikzcd}
		0 \arrow{r} & \RR \arrow[hookrightarrow]{rr} \arrow[hookrightarrow]{d} && \Quad\affi{A} \arrow{rr}{\nabla} \arrow{rrd}{\Hess{}} &&\Aff^S(\affi{A},\vect{\affi{A}}^*) \arrow{d}{\dd} \arrow{r} & 0 \\
		0 \arrow{r}& \Aff\affi{A} \arrow[hookrightarrow]{rru} \arrow[hookrightarrow,swap]{rr}{\Id\otimes u+u\otimes \Id}  && S^2\Aff\affi{A} \arrow[swap]{rr}{\dd\otimes\dd} \arrow{u}{\rotatebox{90}{\(\sim\)}} && S^2\vect{\affi{A}}^* \arrow{r}& 0.
	\end{tikzcd}
	\end{equation}	
	This isomorphism $\Quad\affi{A}\simeq S^2\Aff\affi{A}$ relates quadratic functions $\delta\in\Quad\affi{A}$ with symmetric bilinear forms $\widehat{g}\in S^2\Aff\affi{A}$ on $(\Aff\affi{A})^*=\widehat{\affi{A}}$, and using a referential $\mathcal{R}=(R_0,\ldots,R_n)$, with $(n+1)\times (n+1)$ symmetric matrices $G\in\Sym_{n+1}(\RR)$ (Gram matrix of $\widehat{g}$ in the basis $z_{R_0},\ldots, z_{R_n}$) following the rules:
	\begin{equation}\label{formulagtildezetas}
		\begin{aligned}
			\delta(P)&=\frac12 \widehat{g}(z_P,z_P)=\frac12 p^t\cdot G\cdot p,	\\
			\widehat{g}(z_P,z_Q)&=4\delta\left( \frac{P+Q}{2}\right)-\delta(P)-\delta(Q)=p^t\cdot G\cdot q,\\
			G_{ij}&=\widehat{g}(z_{R_i},z_{R_j})=4\delta\left( \frac{R_i+R_j}{2}\right)-\delta(R_i)-\delta(R_j),
	\end{aligned}\end{equation}
	where $p,q\in\WW_n$ represent the barycentric coordinate vectors associated to points $P,Q$ in the referential $\mathcal{R}$, and $G$ is Gram matrix associated to $\widehat{g}$ in the basis $z_{R_0},\ldots, z_{R_n}$.
\end{theorem}

Theorem \ref{largepropo} shall be proven following next remarks. 
\begin{remark} \label{restringe}
Recall that $\Aff{\affi{A}}=\widehat{\affi{A}}^*$. Hence affine functions on $\affi{A}$ can be seen as the restriction to $\affi{A}\subseteq \widehat{\affi{A}}$ of linear functions defined on $\widehat{\affi{A}}$ (using the immersion $P\mapsto z_P$). The differential $\dd$ is simply the restriction of this linear function to the subspace $\vect{A}\subset\widehat{\affi{A}}$. In the same manner the mapping $\dd\otimes\dd$ that determines the hessian principal component $\Hess{g}$ can be seen as the symmetric bilinear form $\widehat{g}$  associated to $\delta$ restricted to the subspace $\vect{A}=\langle u\rangle^\circ$. The gradient and hessian mappings are also natural. They naturally lead to the factor $1/2$ used in remark \ref{quadraticfrombilinear} to define (together with the immersion $z$) the vertical isomorphism $S^2\Aff{\affi{A}}\rightarrow \Quad\affi{A}$ in (\ref{hessianexactsequence}). Other choices lead to a non-commutative diagram. 

Sticking to exact sequence (\ref{hessianexactsequence}) leads to an identification of the classical bilinear form $\dd x\otimes\dd x+\dd y\otimes \dd y+\dd z\otimes \dd z$ in $\RR^3$ with the equivalence class of quadratic functions $\frac12(x^2+y^2+z^2)+c_1x+c_2y+c_3z+d$, which uses the $1/2$-factor.
 
 Matrix $G$ in (\ref{formulagtildezetas}) is called the linear hull Gram matrix representation of the quadratic function $\delta$, in the affine referential $\mathcal{R}$. It is related to baricentric matrix representation $\Delta$  of $\delta$ in the same referential given in (\ref{matrizDelta}) by $\Delta=\frac12 G$.
\end{remark}
\begin{proof}[Proof of theorem \ref{largepropo}]
The exactness of the upper row in (\ref{hessianexactsequence}) was given by affine Poincare's lemma (proposition \ref{poincare}). 
For the exactness of the lower row: As $\dd\colon \Aff{\affi{A}}\rightarrow \vect{A}^*$ is surjetive, so is $\dd\otimes\dd$. On the other hand, as $u\neq 0$, there holds $f\otimes u+u\otimes f=0$ only when $f=0$. Therefore $\Id\otimes u+u\otimes\Id$ is injective. As $\dd u=0$, the image of $\Id\otimes u+u\otimes \Id$ is contained in the nullspace of $\dd\otimes \dd$. Computing the dimensions we conclude that the lower row is exact.

The commutativity of the upper left hand side triangle is evident, as all morphisms are the natural inclusion morphisms. The commutativity of the upper right hand side triangle is evident, by definition of $\Hess\delta$ as the linear principal component of the affine covector field $\nabla\delta$.

 For the commutativity of the lower left hand side triangle: For any $f\in\Aff{\affi{A}}$ the element $f\otimes u+u\otimes f\in S^2\Aff{\affi{A}}$ determines the quadratic function $\delta(P)=\frac12(f\otimes u+u\otimes f)(z_P,z_P)=\frac12(f(P)\cdot u(P)+u(P)\cdot f(P))=f(P)$, hence the quadratic function is precisely the original affine function $f$.

For the commutativity of the lower right hand side triangle: From $\widehat{g}=(f\otimes h+h\otimes f)\in S^2\Aff{\affi{A}}$ (determined by two affine functions $f,h\in\Aff{\affi{A}}$) we get the quadratic function $\delta(P)=\frac12\widehat{g}(z_P,z_P)=f(P)\cdot h(P)$. Hence $\delta=f\cdot h$. The associated homogeneous quadratic function at $P$ is then $\delta^h_P=(f-f(P))\cdot (h-h(P))$. Hence $\delta-\delta^h_P=f(P)\cdot h+h(P)\cdot f-f(P)\cdot h(P)$. Taking the quotient in $\Aff\affi{A}/\RR=\vect{\affi{A}}^*$ we deduce that the gradient covector field is the affine mapping $P\mapsto \nabla_P\delta=f(P)\cdot \dd h+h(P)\cdot \dd f$, whose principal component is $\Hess\delta=\dd f\otimes\dd h+\dd h\otimes \dd f=(\dd\otimes\dd)(f\otimes h+h\otimes f)$, as we wanted to prove. As these particular bilinear forms $f\otimes h+h\otimes f$ span the whole space $S^2\Aff{\affi{A}}$, we conclude the lower right hand side triangle is commutative.

To prove that $S^2\Aff\affi{A}\rightarrow \Quad\affi{A}$ is injective, observe that any element $\widehat{g}$ that is transformed into $\delta=0$ has null hessian, hence due to the comutativity of our diagram $(\dd\otimes\dd)(\widehat{g})=0$. As the lower row is exact this means that $\widehat{g}=f\otimes u+u\otimes f$ for some affine function $f\in\Aff{\affi{A}}$. But we also know that $f\otimes u+u\otimes f$ determines the quadratic function $\delta(P)=\frac12(f(P)\cdot u(P)+u(P)\cdot f(P))=f(P)$. Hence $\widehat{g}$ determines the null quadratic function $\delta=0$ if and only if $\widehat{g}=0\otimes u+u\otimes 0=0$. The mapping taking $\widehat{g}$ to $\delta$ is linear, injective, and by dimension computation, it is an isomorphism.

We observe further that $P\mapsto z_P$ is an affine transformation, hence $$z_{\frac12(P+Q)}=\frac12(z_P+z_Q),$$ and by symmetry and bilinearity of $\widehat{g}$ we get for the associated $\delta\in\Quad\affi{A}$:
$$	\begin{aligned} 4\delta&\left( \frac{P+Q}{2}\right)-\delta(P)-\delta(Q)=\\ &=
2\widehat{g}\left(\frac12(z_P+z_Q),\frac12(z_P+z_Q)\right)-\frac12\widehat{g}(z_P,z_P)-\frac12\widehat{g}(z_Q,z_Q)=\widehat{g}(z_P,z_Q).\end{aligned}$$

Using this formula, together with the (defining) identities relating $G$ and $\delta$ to $\widehat{g}$, all the rules given in (\ref{formulagtildezetas}) are straightforward.
 \end{proof}
In this formulation, if two affine functions $f,h\in\Aff\affi{A}$ have coefficient row vectors $a,b\in M_{1(n+1)}(\RR)$ with respect to some affine referential, then $f\cdot h$ as a quadratic function has a linear hull Gram matrix representation $G=a^t\cdot b+b^t\cdot a$. In particular the affine function with coefficient row vector $c\in M_{1(n+1)}(\RR)$ has a linear hull Gram matrix representation (as quadratic function)  \begin{equation*} 
	 G_c=c^t\cdot \one_n+\one_n^t\cdot c.\end{equation*} 
 
\section{Metrics on affine spaces}\label{sec4}

Symmetric bilinear forms on a vector space lead to quadratic functions on it that, in particular, turn out to be homogeneous (with respect to the zero vector). Recall that in the affine setting there is no possibility to single out a family of homogeneous quadratic functions. In other words, there is no natural splitting choice for the exact sequence of vector space morphisms:
\begin{equation}\label{SucExMet} 0\rightarrow \Aff{\affi{A}}\hookrightarrow \Quad{\affi{A}}\rightarrow  
\Quad{\affi{A}}/\Aff{\affi{A}}\rightarrow 0.\end{equation}

\begin{define}\label{definemetricas}
	The quotient vector space $\Quad\affi{A}/\Aff\affi{A}$ shall be called space of metrics on the affine space $\affi{A}$, and represented as $\Met\affi{A}$. Its fibers are affine subspaces $\Quad_m\affi{A}$, whose elements are called $m$-quadratic functions (all of them share the same hessian principal component, the hessian principal component $g_m$ associated to the metric $m$).
\end{define}
From exact sequences (\ref{hessianexactsequence}) in theorem \ref{largepropo} we conclude that a metric $m\in\Quad\affi{A}/\Aff\affi{A}$ can be identified using the hessian principal component $g=\Hess\delta$ associated to any of its representatives $\delta$. Each representative $\delta$ can be identified with a symmetric bilinear form $\widehat{g}$ on the linear hull $\widehat{\affi{A}}$.  Following remark \ref{restringe}  the restriction of $\widehat{g}$ to  $\vect{A}\subset\widehat{\affi{A}}$ is $g$.

There are different choices of splitting for the exact sequence (\ref{SucExMet}).
	\begin{itemize}
		\item 	For any choice of point $P\in\affi{A}$, the linear transformation $\delta\in\Quad\affi{A}\rightarrow f\in\Aff\affi{A}$ that transforms $\delta$ into the affine function $\delta-\delta^h_P$ is a retraction of the immersion in the exact sequence (\ref{SucExMet}).
		
		\item 	For any choice of referential $\mathcal{R}=(R_0,\ldots, R_n)$ on $\affi{A}$, the linear transformation $\delta\in\Quad\affi{A}\mapsto f\in\Aff\affi{A}$ that transforms any quadratic function $\delta$ into the unique affine function $f$ such that $f(R_i)=\delta(R_i)$ is a retraction of the immersion in the exact sequence (\ref{SucExMet}).
	\end{itemize}
Considering the difference of any quadratic function $\delta$ with the affine function obtained with any of the retractions above, one gets new quadratic functions $\delta^h_P$ ($P$-homogeneous component of $\delta$), or $\delta^0_{\mathcal{R}}$ ($\mathcal{R}$-reduced component of $\delta$). The space of metrics has then a natural identification with the nullspace of the chosen retraction:
$$\Met\affi{A}\simeq \Quad^h_P\affi{A}=\{\delta\in\Quad\affi{A}\,\colon\, \delta\text{ homogeneous at }P\},$$
$$\Met\affi{A}\simeq \Quad^0_{\mathcal{R}}\affi{A}=\{\delta\in\Quad\affi{A}\,\colon\, \delta\text{ vanishes at }\mathcal{R}\}.$$

Using Gram representation of quadratic functions (theorem \ref{largepropo}) and the specific immersion  $\Met\affi{A}\simeq \Quad^0_{\mathcal{R}}\affi{A}\hookrightarrow  \Quad\affi{A}$ determined by an affine referential $\mathcal{R}$, metrics are identified as symmetric matrices with vanishing diagonal:
$$\Quad^0_{\mathcal{R}}\simeq \{G\in\Sym_{n+1}(\RR)\,\colon\, \diag G=\zero_n\}:=\HH\Sym_{n+1}(\RR)$$
where $\diag\colon M_{kk}(\RR)\rightarrow M_{1k}(\RR)$ represents the identification of the diagonal (as a row vector) and $\zero_n\in M_{1n+1}(\RR)$ is the null row vector.
Symmetric matrices with null diagonal are called hollow symmetric matrices.

The space of metrics is identified with the subspace of quadratic functions such that $\delta=\delta^h_P$. We may consider Gram representation of quadratic functions with respect to some affine referential $\mathcal{R}$ and the specific immersion  $\Met\affi{A}\hookrightarrow \Quad^h_P\affi{A}\subset \Quad\affi{A}$ determined by some point $P\in\affi{A}$ with barycentric coordinate vector $p\in\WW_n$. Using (\ref{coefdeltap}), the relation $G=\frac12\Delta$ given in remark \ref{restringe}, and taking into account that $2p^tG=p^tGp\one$ only happens when $p^tG=\zero_n$ (multiply on the right with $p$, for which $\one p=1$ holds) we get:
\begin{equation}\label{homogeneamatriz}\Quad^h_P\simeq \{G\in\Sym_{n+1}(\RR)\,\colon \, p^tG=\zero_n\}=:\mathbb{Z}_p\Sym_{n+1}(\RR).\end{equation}

Retractions of (\ref{SucExMet}) determine corresponding sections:
\begin{prop}\label{formulametrica}
	Consider any metric $m\in\Met\affi{A}=\Quad\affi{A}/\Aff\affi{A}$.
		\begin{itemize}
			\item For each point $P\in\affi{A}$ there exists a unique $m$-quadratic function $\delta^m_P\in\Quad\affi{A}$ that is homogeneous at $P$ ($P$-homogeneous representative of $m$).
			\item For each referential $\mathcal{R}=(R_0,\ldots, R_n)$ there exists a unique $m$-quadratic function $\delta^m_{\mathcal{R}}\in\Quad\affi{A}$ that vanishes at every point of the referential ($\mathcal{R}$-reduced representative of $m$).
		\end{itemize}  
\end{prop}
\begin{define}\label{def_deltaR_deltam}
	We call half-squared pseudodistance function associated to the metric $m\in\Met\affi{A}$ with respect to the point $P\in\affi{A}$ the only quadratic function $\delta^m_P\in\Quad\affi{A}$ representing the metric and homogeneous with respect to $P$.

	We call $\mathcal{R}$-reduced quadratic function associated to the metric $m\in\Met\affi{A}$ with respect to the referential $\mathcal{R}=(R_0,\ldots, R_n)$ of $\affi{A}$ the only quadratic function $\delta^m_\mathcal{R}\in\Quad\affi{A}$ representing the metric and vanishing at each  point of the referential.
\end{define}
The function $\delta^m_P$ can be obtained taking the homogenization at $P$ of any representative function $\delta\in\Quad\affi{A}$ of the metric $m\in\Quad\affi{A}/\Aff\affi{A}$. As we saw in lemma \ref{lemahesshomogeneo}:
$$\delta^m_P(Q)=\frac12 g(\vect{PQ},\vect{PQ})=\frac12 \widehat g(z_{Q}-z_{P},z_{Q}-z_{P}),$$
using the hessian principal component $g\in S^2\vect{A}^*$ associated to the metric, or the symmetric bilinear form $\widehat{g}\in S^2\Aff\affi{A}$ associated to any representative $\delta\in\Quad\affi{A}$ of the metric $m$.

The function $\delta^m_\mathcal{R}$ can be obtained taking any representative function $\delta\in\Quad\affi{A}$ of the metric $m\in\Quad\affi{A}/\Aff\affi{A}$ and substracting the unique affine function that shares with $\delta$ the same values at every point $R_i$ of the referential.

For any $m$-quadratic function $\delta$, if $G\in\Sym_{n+1}(\RR)$ is Gram matrix associated to $\delta$ in the referential $\mathcal{R}=(R_0,\ldots,R_n)$ we may use the relation $G=2\Delta$ with the barycentric matrix representation of $\delta$ and formula (\ref{formula4b}) to prove that $\delta^m_P$ has, in the same referential, Gram matrix:
	\begin{equation}\label{HalfsquaredfromGram}
		G^h_p=(\Id_n-p\cdot \one_n)^t\cdot G\cdot (\Id_n-p\cdot \one_n)\in \mathbb{Z}_p\Sym_{n+1}(\RR),
	\end{equation}
where $p\in\WW_n$ is the barycentric coordinate vector associated to $P$ in the referential $\mathcal{R}$.	Using (\ref{coefdelta0})  we get Gram matrix associated to $\delta^m_{\mathcal{R}}$: 	$$G^0=G-\frac12\left(\one_n^t\cdot \diag G+(\diag G)^t\cdot \one_n\right)\in\HH\Sym_{n+1}(\RR).$$

\begin{remark}
	An affine function $f$ is in particular an element $f\in\Quad\affi{A}$, associated by theorem \ref{largepropo} with a symmetric bilinear form $\widehat{g}=f\otimes u+u\otimes f$:
	$$\begin{aligned} \widehat{g}(z_{R_i},z_{R_j})&=(f\otimes u+u\otimes f)(z_{R_i},z_{R_j})=f(R_i)+f(R_j).
	\end{aligned}$$

	Hence if we take any $m$-quadratic function $\delta\in\Quad_m\affi{A}$ and $f$ affine such that $f(R_i)=\delta(R_i)$, we may represent the unique quadratic function $\delta^m_{\mathcal{R}}=\delta-f$ that is metric-equivalent with $\delta$ and vanishes on the referential $\mathcal{R}$. It determines following (\ref{formulagtildezetas}) a symmetric bilinear form $\widehat{g}_{\mathcal{R}}$ on $\widehat{\affi{A}}$, using $z_{R_0},\ldots, z_{R_n}$ as basis for $\widehat{\affi{A}}$. We get:
	\begin{equation}\label{g0c0}
		\begin{aligned}
	\widehat{g}_{\mathcal{R}}^m(z_{R_i},z_{R_j})&=\left(4\delta\left( \frac{R_i+R_j}{2}\right)-\delta(R_i)-\delta(R_j)\right)-\left(\delta(R_i)+\delta(R_j)\right)=\\&=4\delta\left( \frac{R_i+R_j}{2}\right)-2\delta(R_i)-2\delta(R_j),\end{aligned}
	\end{equation}
	which is represented by a Gram matrix with null diagonal entries.
\end{remark}

Any affine mapping $\affi{B}\rightarrow \affi{A}$ induces $\Quad\affi{A}\rightarrow \Quad\affi{B}$ that transforms quadratic functions into quadratic functions, taking the subspace of affine functions into the subspace of affine functions. Therefore it determines a linear mapping $\Met\affi{A}\rightarrow \Met\affi{B}$ on the corresponding quotient spaces. Hence, given any affine immersion $\affi{B}\subseteq \affi{A}$, we may restrict any metric $m\in\Met\affi{A}$ to a metric $\left. m\right|_B\in\Met\affi{B}$, or given any affine projection $\affi{B}\rightarrow \affi{A}$ we may pull-back any metric $m\in\Met\affi{A}$ to a larger affine space $\affi{B}$.
\begin{define}
	The radical of a metric $m\in\Met{\affi{A}}$ is the maximal vector subspace $\vect{R}\subseteq\vect{A}$ such that $m$ is the pull-back of a metric on $\affi{A}/\vect{R}$. A non-degenerate metric is one with a trivial radical.
\end{define}
If $m$ is determined by the self-adjoint bilinear form $g\colon\vect{A}\rightarrow\vect{A}^*$, then the radical of $m$ is the null-space of the bilinear form, the subspace of vectors $\vect{x}\in\vect{A}$ that are transformed into the null covector by the bilinear form.

\begin{prop} \label{proposition_immersionCayley}
	For a given metric $m$ on $\affi{A}$ determined by a symmetric bilinear form $g$, there holds $\delta^m_P=\delta^m_Q$ if and only if $\vect{PQ}$ is in the nullspace of the bilinear form $g$.
\end{prop}
\begin{proof}
	Observe that:
	$$\delta^m_Q(R)=\frac12g(\vect{QR},\vect{QR}),$$
	$$\begin{aligned} \delta^m_P(R)&=\frac12g(\vect{PR},\vect{PR})=\frac12g(\vect{PQ}+\vect{QR},\vect{PQ}+\vect{QR})=\\ &= \frac12 g(\vect{PQ},\vect{PQ})+g(\vect{PQ},\vect{QR})+\frac12g(\vect{QR},\vect{QR}).\end{aligned}$$

	Hence stating $\delta^m_P=\delta^m_Q$ is equivalent to state that:
	$$g(\vect{PQ},\vect{QR})=-\frac12g(\vect{PQ},\vect{PQ}),\quad \forall \vect{QR}\in\vect{A}.$$
	
	By linearity, the only way that $g(\vect{PQ},\vect{QR})$ could be independent of $R$ is that $\vect{PQ}$ is in the nullspace of $g$, as we wanted to prove.
\end{proof}
\begin{corollary}\label{coroCMimmersion}
	The mapping $P\in\affi{A}\mapsto \delta^m_P\in\Quad_m\affi{A}$ is a quadratic mapping. It is an immersion if an only if $m$ has trivial radical. 
\end{corollary}
\begin{proof}
	Observe that $\delta^m_X(Y)=\delta^m_Y(X)$	for any points $X,Y\in\affi{A}$. Let us prove that the mapping $P\mapsto \delta^m_P$ is quadratic. For any pair of points $P,Q$ and $1$-dimensional weight $(\alpha,\beta)$ the image of $\alpha\cdot P+\beta\cdot Q$ is the quadratic function:
	$$\delta^m_{\alpha P+\beta Q}(X)=\delta^m_X(\alpha P+\beta Q)=\ldots$$
	As $\delta^m_X$ is quadratic we get:
	$$\ldots=\alpha(\alpha-\beta)\delta^m_X(P)+\beta(\beta-\alpha)\delta^m_X(Q)+4\alpha\beta \delta^m_X\left(\frac{P+Q}{2}\right).$$
	Hence observing again that $\delta^m_X(Y)=\delta^m_Y(X)$ we conclude:
	$$\delta^m_{\alpha P+\beta Q}(X)=\left(\alpha(\alpha-\beta)\delta^m_P+\beta(\beta-\alpha)\delta^m_Q+4\alpha\beta \delta^m_{\frac{P+Q}{2}}\right)(X),$$
	which represents the quadraticity condition for the mapping $\delta^m\colon\affi{A}\rightarrow \Quad_m\affi{A}$.
	
	The equivalence of the injectivity of the mapping $\delta^m$ and non-degeneracy of the metric $m$ was proven in proposition \ref{proposition_immersionCayley}
\end{proof}
\begin{define}\label{coroCMimmersiondef}
	For a non-dentenerate metric $m\in\Met\affi{A}$ we call $\delta^m\colon\affi{A}\rightarrow \Quad\affi{A}$ Cayley-Menger immersion of the affine space $\affi{A}$.
\end{define}
\begin{remark}
We may identify the coordinate representation of the quadratic immersion $P\mapsto \delta^m_P$. Fix an affine referential $\mathcal{R}=(R_0,\ldots, R_n)$ on $\affi{A}$ and consider the hollow symmetric matrix $G$ representing $m$ in this referential. Each element $P\in\affi{A}$ is determined by its barycentric coordinate vector $p=[p_0\ldots p_n]^t\in\WW_n$. Formula (\ref{HalfsquaredfromGram}) shows that $\delta^m_P$ has dependence on $P$ (corollary \ref{coroCMimmersion}) and linear dependence on $m$.
\end{remark}
Recall that the notion of convexity of a quadratic function (definition \ref{defconvex}) is preserved when an affine function is added. Therefore it is a notion that can be defined for equivalence classes $m\in \Quad\affi{A}/\Aff\affi{A}$.
\begin{define}
We say the metric $m\in\Met\affi{A}=\Quad\affi{A}/\Aff\affi{A}$ is positive definite if some/any of its representatives $\delta\in\Quad\affi{A}$ is convex. We say a metric is null if some/any of its representatives is an affine function.
\end{define}

If $g_m$ is the hessian principal component of a metric $m$, then positivity/nullity of the metric on an affine subspace $F$ corresponds, respectively, to positive-definiteness and vanishing of the bilinear form $g_m$ restricted to the director vector subspace $\vect{F}$. Affine spaces where a metric has null restriction are called isotropic.
\begin{define}
	The dimension of a maximal affine subspace where $m$ restricts as positive definite is called positivity $\pi(m)$ of the metric. The dimension of a maximal affine subspace where $-m$ restricts as positive definite is called negativity $\nu(m)$ of the metric. The dimension of the radical is called nullity $\rho(m)$ of the metric. The integer values $(\pi(m),\nu(m),\rho(m))$ form the so-called inertia index of the metric.
\end{define}

\begin{remark}\label{remark43}
	Consider any metric $m\in\Met\affi{A}$ on the affine-space $\affi{A}$, with associated hessian principal component $g\in S^2\vect{A}^*$. Using exact sequence (\ref{hessianexactsequence}) all quadratic functions $\delta$ representing the metric can be identified as elements $\widehat{g}\in S^2\Aff\affi{A}=S^2\widehat{\affi{A}}^*$, and the restriction to $\vect{A}\subset\widehat{\affi{A}}$ is the symmetric bilinear form $g$ on $\vect{A}$ associated to $m$.
	
	The study of the inertia of the metric $m\in\Met\affi{A}$ is equivalent to the study of the inertia of the restriction to $\langle u\rangle^\circ=\vect{A}$ of the symmetric bilinear form $\widehat{g}$ on $\widehat{\affi{A}}$ associated to any of its representatives $\delta$. 
	The main tools to study inertia of bilinear forms restricted to subspaces can be found in \cite{maddocks1988restricted}. For quadratic functions that are homogeneous at a point $P$ the relation is straightforward: Following (\ref{homogeneamatriz}), the bilinear form $\widehat{g}$ on $\widehat{A}$ represents a quadratic function that is homogeneous at $P\in\affi{A}$ precisely when $z_P$ is in its nullspace. As $\vect{A}\oplus\langle z_P\rangle=\widehat{\affi{A}}$ and $\widehat{g}$ restricts as the hessian principal component $g$ on the subspace $\vect{A}$, for $g$ with inertia index $(\pi,\nu,\rho)$ and homogeneous at $P$ we have $\widehat{g}$ with inertia index $(\pi,\nu,\rho+1)$.
		
	If we choose an affine referential $\mathcal{R}=(R_0,\ldots, R_n)$, we have a basis $z_{R_0},$ $\ldots,$ $z_{R_n}$ of $\widehat{\affi{A}}$ and formula (\ref{formulagtildezetas}) determines the Gram matrix $G_{ij}=4\delta(R_{ij})-\delta(R_{i})-\delta(R_{j})$ associated to the bilinear quadratic form $\widehat{g}$ corresponding to $\delta\in\Quad_m\affi{A}$. The element $u$, which is a linear function on $\widehat{\affi{A}}$ is represented by the row vector $\one_n$.
	
	A necessary and sufficient condition for $\delta\in\Quad\affi{A}$ to represent a non-degenerate metric is that the following symmetric matrix is non-degenerate:
	$$\left[ \begin{matrix} -G&\one_n^t \\ \one_n & 0\end{matrix}\right],\qquad G_{ij}=4\delta\left(\frac{R_i+R_j}{2}\right)-\delta(R_i)-\delta(R_j).$$
	
	Indeed, for $g$ to be degenerate, there should exist an element in $\langle\one_n\rangle^\circ$ that is $G$-orthogonal to every other element of $\langle\one_n\rangle^\circ$, therefore there should exist a non-trivial solution to the system of equations $\one_n\cdot x=0$, $G\cdot x=\alpha\cdot \one_n^t,$
	which can be represented as the matrix equation:
	$$\left[\begin{matrix} -G&\one_n^t\\ \one_n & 0\end{matrix}\right]\cdot \left[\begin{matrix} x\\ \alpha\end{matrix}\right]=0.$$
	
	Nontrivial solutions of this matrix equation determine the radical of the metric. 

	The matrix $G$ above depends on a specific choice of quadratic function $\delta$ that represents the metric $m$, together with a choice of affine referential $\mathcal{R}$. We know that in fact this choice of referential leads to a specific quadratic function $\delta^m_{\mathcal{R}}$ representing the metric $m$. This function vanishes at all the referential points. Recall from (\ref{g0c0}) that using any $\delta\in\Quad_m\affi{A}$ its $\mathcal{R}$-reduced component  $\delta^0_{\mathcal{R}}=\delta^m_{\mathcal{R}}$ has an associated bilinear form $\widehat{g}^m_{\mathcal{R}}$ on $\widehat{\affi{A}}$ with Gram matrix:
	$$\widehat{g}^m_{\mathcal{R}}(z_{R_i},z_{R_j})=G^0_{ij}=4\delta\left(\frac{R_i+R_j}{2}\right)-2\delta(R_i)-2\delta(R_j).$$
	
	As $g$ is the restriction of $\widehat{g}^m_{\mathcal{R}}$ to the subspace $\vect{A}$, and as $\widehat{g}^m_{\mathcal{R}}(z_{R_i},z_{R_i})=0$, we get:
	\begin{equation}\label{formulaG} G^0_{ij}=-\frac12 \widehat{g}(z_{R_j}-z_{R_i},z_{R_j}-z_{R_i})=-\frac12 g(\vect{R_iR_j},\vect{R_iR_j}).\end{equation}
\end{remark}
\begin{prop}\label{nondegeneracy}
Fix a referential $\mathcal{R}=(R_0,\ldots,R_n)$ on the affine space $\affi{A}$. A metric $m\in\Met\affi{A}$ is non-degenerate if and only if the following hollow symmetric matrix is non-degenerate:
	\begin{equation}\label{CMmatrix}
		\left[\begin{matrix} D_{ij} & \one_n^t \\ \one_n & 0\end{matrix}\right]\qquad D_{ij}=g(\vect{R_iR_j},\vect{R_iR_j}).
	\end{equation}
\end{prop}
\begin{proof}
It is a simple application of the previous remark, using $\delta^m_{\mathcal{R}}$ as quadratic function representing the metric, and multiplying first rows with $-1/2$ and last column with $-2$.
\end{proof}
Observe that for any metric $m$ and referential $\mathcal{R}$ with  squared distance matrix $D$, and for the barycentric representation $\Delta$ and Gram representation $G$ associated to $\delta^m_{\mathcal{R}}$ there holds $2\Delta=G=-\frac12 D$.

\begin{remark}\label{remarkEDM}
	The first appearance of matrix (\ref{CMmatrix}) goes back to certain results by Cayley \cite{cayley1841theorem}, who studied its degeneracy for any choice of $n+3$ points on $n$-dimensional euclidean spaces. This matrix appears later in Menger's work \cite{menger1930untersuchungen} where several results are proven regarding the isometric immersion of points in euclidean spaces. In euclidean geometry the matrix is classically known as Cayley-Menger matrix associated to the euclidean metric $m$ and to the referential $\mathcal{R}=(R_0,\ldots, R_n)$. Several other works deal with this determinant and matrix. For the case of positive definite metrics, the matrix $D$ is called Euclidean Distance Matrix (EDM). The interested reader may consult \cite{alfakih2018euclidean} on this subject. For a treatment of the Cayley-Menger determinant (and bi-determinant) from an euclidean geometry perspective, the reader is referred to section 9.7 in Berger's Geometry notes \cite{berger2009geometry}
	
	All these works deal with the EDM or the Cayley-Menger matrix from the perspective of matrix calculus. We have in (\ref{formulaG}) an additional meaning of the EDM (up to a factor $-1/2$) as Gram matrix associated to a specific symmetric bilinear form on $\widehat{\affi{A}}$ (determined by $m$ and a choice of referential). It is the Gram matrix representation of the unique $m$-quadratic function $\delta^m_{\mathcal{R}}$ that vanishes at the referential. The symmetric EDM matrix $D$ depends on the choice of referential and is not associated to any intrinsically defined metric on $\widehat{\affi{A}}$.
 \end{remark}
		
\section{The Cayley-Menger bilinear form associated to a metric}\label{sec5}
Recall that each affine space $\affi{A}$ has an associated vector space $\widehat{\affi{A}}$ (the linear hull) and an affine immersion $z\colon\affi{A}\hookrightarrow \widehat{\affi{A}}$, so that all affine functions on $\affi{A}$ are obtained as composition of some linear function $f\colon\widehat{\affi{A}}\rightarrow \RR$ with $z$.

In the same manner we may  introduce a quadratic hull associated to $\affi{A}$.
\begin{define}\label{define_vP}
	The vector space $\widehat{\affi{A}}^2=(\Quad\affi{A})^*$ of linear operators on the space of quadratic functions of $\affi{A}$ is called quadratic hull associated to $\affi{A}$. Each point $P\in\affi{A}$ induces an element $v_P$ on $\widehat{\affi{A}}^2$, determined by $\langle v_P,\delta\rangle=\delta(P)$, for any quadratic function $\delta\in\Quad\affi{A}$. The element $v_P$ is called evaluation operator at $P$ on quadratic functions. 
\end{define}
\begin{prop}\label{v_P_isquadratic}
	The mapping $P\in\affi{A}\mapsto v_P\in(\Quad\affi{A})^*=\widehat{\affi{A}}^2$ is a quadratic immersion. The linear inclusion $\Aff\affi{A}\subset \Quad\affi{A}$ determines a linear projector $\widehat{\affi{A}}^2=(\Quad\affi{A})^*\rightarrow (\Aff\affi{A})^*=\widehat{\affi{A}}$, which together with the affine immersion $P\mapsto z_P$ leads to a commutative triangle 
	$$	\begin{tikzcd}
		& \widehat{\affi{A}}^2\arrow{dr}{lin.proj.}& \\
		\affi{A} \arrow[Rightarrow]{ur}{v\,(quad.)} \arrow[hookrightarrow]{rr}{z\, (aff.)} && \widehat{\affi{A}}
	\end{tikzcd}$$
(where double arrows are used to represent quadratic mappings).\end{prop}
\begin{proof}
	To see that $P\mapsto v_P$ is quadratic we use the quadraticity condition (\ref{quadraticity}). For any pair of points $P,Q\in\affi{A}$ and 1-dimensional weight $(\alpha,\beta)\in \WW_1$, and for any quadratic function $\delta\in\Quad\affi{A}$ we observe:
	$$\begin{aligned} \langle v_{\alpha P+\beta Q},\delta\rangle&=\delta(\alpha P+\beta Q)=\\
		&=\alpha(\beta-\alpha)\cdot \delta(P)+\beta(\alpha-\beta)\cdot\delta(Q)+4\alpha\beta\cdot \delta\left(\frac{P+Q}{2}\right)=\\
		&=
		\left\langle \alpha(\beta-\alpha)\cdot v_P+\beta(\alpha-\beta)\cdot v_Q+4\alpha\beta\cdot v_{\frac{P+Q}{2}},\delta\right\rangle.\end{aligned}$$

	Hence $v_{\alpha P+\beta Q}=\alpha(\beta-\alpha)\cdot v_P+\beta(\alpha-\beta)\cdot v_Q+4\alpha\beta\cdot v_{\frac{P+Q}{2}}$ and the mapping $v$ is quadratic.
	
	It is evident that for $f\in\Aff\affi{A}$, there holds $\langle z_P,f\rangle=f(P)=\langle v_P,i(f)\rangle$. Here $i$ represents the inclusion of $\Aff\affi{A}$ into $\Quad\affi{A}$, that induces the linear projector $(\Quad\affi{A})^*\rightarrow (\Aff{\affi{A}})^*$, giving the commutative triangle in our statement.
	
	As $z$ is injective and factors by $v$, we conclude that $v$ is also injective, hence a quadratic immersion.
\end{proof}
Each linear function $\widehat{\affi{A}}^2\rightarrow \RR$, by composition with the quadratic immersion $v\colon \affi{A}\hookrightarrow \widehat{\affi{A}}^2$,  determines a quadratic function on $\affi{A}$. Observe that conversely, any quadratic function $\delta\in\Quad\affi{A}$ represents a linear form on the dual space $(\Quad\affi{A})^*=\widehat{\affi{A}}^2$, that is, a linear mapping $\widehat{\affi{A}}^2\rightarrow \RR$ that transforms any linear functional $F$ defined on $\Quad\affi{A}$ into the real value $\langle F,\delta\rangle$. Moreover, the composition of this linear mapping with the quadratic immersion $v$ is the quadratic function $\delta$ itself: $P\mapsto \langle v_P,\delta\rangle =\delta(P)$.

All quadratic functions are obtained by composition of linear forms of the quadratic hull with the natural quadratic immersion $v\colon \affi{A}\hookrightarrow \widehat{\affi{A}}^2$.

\begin{remark}
	We may use a referential $\mathcal{R}=(R_0,\ldots, R_n)$ to identify $\affi{A}$ with the affine subspace $\WW_n$ taking any point $P\in\affi{A}$ to its barycentric coordinate vector $p=[p_0\ldots p_n]^t$, such that $\one_n\cdot p=1$. As already observed in (\ref{formulagtildezetas}) the same referential allows to identify quadratic functions $\delta$ on $\affi{A}$ with Gram symmetric matrices $G\in\Sym_{n+1}(\RR)$, so that $\delta(P)=\frac12 p^t\cdot G\cdot p$.
	
	If we use the trace-duality $\langle C,D\rangle=\trz(C\cdot D)=\sum c_{ij}\cdot d_{ij}$, the dual space of $\Sym_{n+1}(\RR)$ is identified with itself. In this case, observing that
	$$\langle \delta,v_P\rangle=\frac12p^t\cdot G\cdot p=\frac12\sum G_{ij}p_i\cdot p_j=\sum G_{ij}\cdot H_{ij},$$
	it suffices to take the matrix with entries $H_{ij}=\frac12p_i\cdot p_j$. Hence in our o.d.system of affine coordinates the quadratic immersion $P\mapsto v_P$ is written as $p\mapsto \frac12 p\cdot p^t$, clearly a quadratic expression on the affine space $\affi{A}$.
\end{remark}

Consider a metric $m\in\Quad\affi{A}/\Aff\affi{A}$ on some $n$-dimensional affine space $\affi{A}$. The fiber $\Quad_m\affi{A}\subset \Quad\affi{A}$ is then an affine subspace, whose director vector subspace is the $n+1$-dimensional subspace $\Aff\affi{A}\subset\Quad\affi{A}$. In other words, the projection $\Quad\affi{A}\rightarrow \Met\affi{A}$ can be seen as an affine bundle on the (vector) space of metrics, whose director vector bundle is the trivial bundle on $\Met\affi{A}$ with fiber $\Aff\affi{A}$.

There exists a linear mapping $\pi_m\colon (\Quad\affi{A})^*\rightarrow \Aff\Quad_m\affi{A}$ taking linear forms (``functionals'') on $\Quad\affi{A}$ to its restriction to the affine subspace $\Quad_m\affi{A}$, which is an affine functional on the affine space of $m$-quadratic functions $\Quad_m\affi{A}$. Therefore there exists a natural quadratic mapping $\pi_m\circ v\colon \affi{A}\rightarrow (\Quad\affi{A})^*\rightarrow \Aff\Quad_m\affi{A}$ that we will still denote by $v$.

Consider now, for any fixed metric $m\in\Met\affi{A}$ the $(n+2)$-dimensional vector space $\Aff\Quad_m\affi{A}$ of affine functions on the $n+1$-dimensional affine space  $\Quad_m\affi{A}$. We shall talk of (affine) functionals on the affine space of $m$-quadratic functions. 
\begin{define}\label{defvP}
Specific elements of $\Aff\Quad_m\affi{A}$  are  $v_P\in\Aff\Quad_m\affi{A}$, the operator defined for any point $P\in\affi{A}$ as restriction of $v_P\in(\Quad \affi{A})^*$ given in definition \ref{define_vP} to the affine subspace $\Quad_m\affi{A}\subset \Quad\affi{A}$, and also the unitary operator $v_m$, defined as $v_m(\delta)=1$ for each $\delta\in\Quad_m\affi{A}$.
\end{define}
\begin{remark}\label{metricanula}
In the specific case of a null metric $m=0\in\Met\affi{A}$, the affine space $\Quad_m\affi{A}$ is identified with $\Aff\affi{A}$ (which is a vector space). In this space we have the zero function $\delta=0\in\Quad_m\affi{A}$. However this function is not the zero element of the linear hull $\widehat{\Quad_m\affi{A}}$, as it belongs to the affine subspace $\Quad_m\affi{A}\subset \widehat{\Quad_m\affi{A}}$ and there holds $\langle \delta, v_m\rangle=1$ in this case. This inconvenience appears always when we consider the linear hull of an affine space that has a vector space structure. In such situations we must distinguish between the zero element of the linear hull and the point determined by the zero vector of the original vector (and affine) space.
\end{remark}
The gradient covector mapping $\dd_B\colon \Aff\affi{B}\rightarrow \vect{B}^*$ in the particular case $\affi{B}=\Quad_m\affi{A}$ (where $\vect{B}=\Aff\affi{A}$) determines a natural linear projector and immersion:
	\begin{equation}\label{define_d_dt}
\dd_m\colon \Aff\Quad_m\affi{A}\rightarrow (\Aff\affi{A})^*,\qquad \dd_m^*\colon\Aff\affi{A}\hookrightarrow \widehat{\Quad_m\affi{A}},
	\end{equation}
where $\widehat{\Quad_m\affi{A}}=\left(\Aff\Quad_m\affi{A} \right)^*$ is the linear hull   (definition \ref{definelinearhull}) of the affine space $\Quad_m\affi{A}$.

Recall that the nullspace of the projector $\dd_m$ is the subspace of constant functionals $\RR=\langle v_m\rangle$, and therefore the image of $\dd_m^*$ is the space of linear forms on $\Aff\Quad_m\affi{A}$ that vanish on  $v_m\in \Aff\Quad_m\affi{A}$.

More specifically, for any affine functional $\bar{v}$ on $\Quad_m\affi{A}$ and any affine function $s\in\Aff\affi{A}$, there holds:
\begin{equation} \label{d_on_functionals} \langle \dd_m \bar{v},s\rangle=\langle \dd_m^*s,\bar{v}\rangle= \bar{v}(\delta+s)-\bar{v}(\delta),\qquad \text{ any }\delta\in\Quad_m\affi{A}. \end{equation}

Recall that $\Quad_m\affi{A}$ can be seen as an affine hyperplane on its linear hull $\widehat{\Quad_m\affi{A}}$. On this linear hull there is a particular element:
\begin{define}
The unit affine function $u\in\Aff\affi{A}$ determines an element $u_m=\dd_m^* u=u\circ\dd_m$ that we call unit element of the linear hull $\widehat{\Quad_m\affi{A}}$.
	
	Affine functionals $\bar{v}\in\Aff\Quad_m\affi{A}$ such that $\langle u_m,\bar{v}\rangle=1$ shall be called normalized affine functionals. They form an affine subspace of functionals:
	$$\bar{v}\in\Aff_1\Quad_m\affi{A}\Leftrightarrow \bar{v}(\delta+\alpha\cdot u)=\alpha+\bar{v}(\delta),\,\forall\delta\in\Quad_m\affi{A},\, \alpha\in\RR.$$
	
	Affine functionals $\bar{v}\in\Aff\Quad_m\affi{A}$ such that $\langle u_m,\bar{v}\rangle=0$ are functionals that factor by the affine quotient space $\Quad_m\affi{A}/\RR$, forming a subspace $\Aff_0\Quad_m\affi{A}$ that shall be called space of $\RR$-projectable affine functionals. 
\end{define}
	The space $\Aff_0\Quad_m\affi{A}$ of affine functionals on the affine space $\Quad_m\affi{A}/\RR$  of $m$-covector fields is the director space of the affine space of normalized affine functionals:
$$\bar{v}\in\Aff_0\Quad_m\affi{A}=\langle u_m\rangle^\circ\Leftrightarrow \bar{v}(\delta+\alpha\cdot u)=\bar{v}(\delta),\,\forall\delta\in\Quad_m\affi{A},\, \alpha\in\RR.$$
	This is the director space of the affine space $\Aff_1\Quad_m\affi{A}$.
	
	From exact sequence (\ref{hessianexactsequence}), elements in $\Quad_m\affi{A}/\RR$ can be seen as affine covector fields on $\affi{A}$ that project by $\dd$ to the symmetric bilinear form on $\vect{A}$ determined by the  metric $m\in\Met\affi{A}=\Quad\affi{A}/\Aff\affi{A}$. These covector fields are closed (symmetric) because the corresponding hessian principal component is associated to $m$. We may then call $\Quad_m\affi{A}/\RR$ the space of $m$-covector fields on $\affi{A}$.
	
For the constant unitary operator $v_m\in\Aff\Quad_m\affi{A}$ and for $u_m=\dd_m^*u$ there holds
\begin{equation}\label{umvm0}\langle v_m,u_m\rangle=\langle \dd_m v_m,u\rangle=v_m(\delta+u)-v_m(\delta)=1-1=0,\quad \text{any }\delta\in\Quad_m\affi{A},\end{equation}
The constant unitary operator $v_m$ belongs to $\Aff_0\Quad_m\affi{A}$, the director vector subspace.

On the other hand for points $P\in\affi{A}$ the functional $v_P\in\Aff\Quad_m\affi{A}$ has the property (any $\delta\in\Quad_m\affi{A}$):
\begin{equation}\label{umvp1}\langle v_P,u_m\rangle=\langle \dd_m v_P,u\rangle=v_P(\delta+u)-v_P(\delta)=\delta(P)+u(P)-\delta(P)=1.\end{equation}
Hence the image of the quadratic immersion $v\colon P\in \affi{A}\hookrightarrow v_P\in\Aff\Quad_m\affi{A}$ is contained in the affine hyperplane $\Aff_1\Quad_m\affi{A}$ of normalized affine functionals.

\begin{lemma}\label{lemadeltav}
	Consider a metric $m\in\Met\affi{A}$ with hessian principal component $g_m\in S^2\vect{A}^*$.  The mapping $\widehat\delta\mapsto \widehat\delta\circ v$ obtained by composition with the quadratic mapping $v\colon \affi{A}\rightarrow \Aff\Quad_m\affi{A}$ is a linear mapping $\widehat{\Quad_m\affi{A}}\rightarrow \Quad\affi{A}$. 
	The restriction of this linear mapping to the hyperplane $\Quad_m\affi{A}\subset \widehat{\Quad_m\affi{A}}$ is the natural immersion $\Quad_m\affi{A}\hookrightarrow \Quad\affi{A}$, and the hessian principal component of $\widehat{\delta}\circ v\in\Quad\affi{A}$ is precisely  $\langle \widehat{\delta},v_m\rangle \cdot g_m$, for any $\widehat{\delta}\in\widehat{\Quad_m\affi{A}}$
\end{lemma}
\begin{proof}
	Any element 
	$\widehat{\delta}\in \widehat{\Quad_m\affi{A}}$ is a linear form on $\Aff\Quad_m\affi{A}$. The composition of this linear form with the quadratic mapping $v\colon \affi{A}\rightarrow \Aff\Quad_m\affi{A}$ is therefore a quadratic mapping from $\affi{A}$ to $\RR$, hence an element $\widehat{\delta}\circ v\in\Quad \affi{A}$. Its hessian principal component is not necessarily $g_m$.
	
	Is is immediate that $(\widehat{\delta}_1+\alpha\widehat{\delta}_2)\circ v=\widehat{\delta}_1\circ v+\alpha\widehat{\delta}_2\circ v$ hence 
	the mapping $\widehat{\delta}\mapsto \widehat{\delta}\circ v$ is a linear mapping.
	
	For any affine space $\affi{B}$ elements $y\in\affi{B}$ determine $\widehat{y}\in\widehat{\affi{B}}=(\Aff\affi{B})^*$ characterized by $\langle \widehat{y},f\rangle=f(y)$. For the particular case $\affi{B}=\Quad_m\affi{A}$ each element $\delta\in\Quad_m\affi{A}$ determines a linear form $\widehat{\delta}\in\widehat{\Quad_m\affi{A}}$ such that $\widehat{\delta}\circ v(P)=\langle \widehat{\delta},v_P\rangle =v_P(\delta)=\delta(P)$. Therefore the mapping $\widehat{\delta}\mapsto \widehat{\delta}\circ v$ restricted to $\Quad_m\affi{A}\subset \widehat{\Quad_m\affi{A}}$ is the identity mapping. The hessian principal component of $\widehat{\delta}\circ v$ is the hessian principal component of $\delta\in\Quad_m\affi{A}$, which is $g_m$, and coincides with $\langle \widehat{\delta},v_m\rangle \cdot g_m$, because $\Quad_m\affi{A}$ is the hyperplane defined by $v_m=1$.
	
	An element that does not belong to the hyperplane $\Quad_m\affi{A}$ is $u_m\in\widehat{\Quad_m\affi{A}}$. Following (\ref{umvp1}) there holds $u_m\circ v(P)=1$. Hence $u_m$ is transformed into the constant unit function $u\in\Quad\affi{A}$. The hessian principal component of this constant function is $0$, which again coincides with $\langle \widehat{\delta},v_m\rangle \cdot g_m$ for the case $\widehat{\delta}=u_m$ (because of (\ref{umvm0})).
	
	Taking the hessian principal component of a quadratic function is a linear mapping. Therefore, the formula proposed for the hessian principal component holds on $\Quad_m\affi{A}\oplus \langle u_m\rangle=\widehat{\Quad_m\affi{A}}$.
\end{proof}
\begin{remark}
For non-vanishing metric $m\in\Met\affi{A}$, the previous lemma determines a linear immersion of $\widehat{\Quad_m\affi{A}}$ into $\Quad\affi{A}$ as the subspace of quadratic functions whose hessian principal component is a multiple of $g_m$. In the case of the null metric $m=0$, this identification is not valid.
\end{remark}

\begin{lemma}\label{projetor_pim}
	The projector $\dd_m\colon \Aff\Quad_m\affi{A}\rightarrow \widehat{\affi{A}}$ described in (\ref{define_d_dt}) transforms the evaluation operator at $P$ on quadratic functions $v_P$ into the evaluation operator at $P$ on affine functions $z_P$. It also transforms the constant unitary operator $v_m$ into $0$.
\end{lemma}
\begin{proof}
	Clearly from (\ref{d_on_functionals}), $$\langle \dd_m v_P,s\rangle=v_P(\delta+s)-v_P(\delta)=(\delta+s)(P)-\delta(P)=s(P)=\langle z_P,s\rangle,\quad \forall s\in\Aff\affi{A}.$$ Hence $\dd_m v_P=z_P\in(\Aff\affi{A})^*=\widehat{\affi{A}}$.
	
	Regarding $v_m$, there holds $\langle \dd_m v_m,s\rangle =v_m(\delta+s)-v_m(\delta)=1-1=0$ for $s\in\Aff\affi{A}$.
\end{proof}
Hence all elements $v_P\in\Aff\Quad_m\affi{A}$ belong to the affine subspace $\dd_m^{-1}(\affi{A})$ (using the immersion $z\colon\affi{A}\hookrightarrow \widehat{\affi{A}}$). Moreover there holds $v_P=v_Q$ only when $z_P=z_Q$, hence for $P=Q$. For any metric $m\in\Met\affi{A}$ the affine space $\affi{A}$ has a natural quadratic immersion into a certain affine hyperplane $\dd_m^{-1}(\affi{A})\subset \Aff\Quad_m\affi{A}$.
\begin{prop}\label{basevR}
	Let $m\in\Met\affi{A}$ be a metric on the affine space $\affi{A}$ and $\mathcal{R}=(R_0,\ldots, R_n)$ any affine referential on this space. The functionals $v_{R_0},$  $\ldots,$ $v_{R_n},$ $v_m$ (see definition \ref{defvP}) form a basis of $\Aff\Quad_m\affi{A}$.

	The corresponding dual basis is $w_0,\ldots, w_n,\delta^m_{\mathcal{R}} \in \widehat{\Quad_m\affi{A}}$ where $w_0,$ $\ldots,$  $w_n\in\Aff\affi{A}$ is the barycentric affine coordinate system associated to the referential, and  $\delta^m_{\mathcal{R}}$ is the $\mathcal{R}$-reduced quadratic function associated to the metric (see definition \ref{def_deltaR_deltam}).
\end{prop}
\begin{proof}
	Consider that a given linear combination of these affine functionals vanishes:
	$$\alpha_0\cdot  v_{R_0}+\ldots+\alpha_n\cdot v_{R_n}+\beta\cdot v_m=0.$$
	
	Consider the projection $\dd_m\colon \Aff\Quad_m\affi{A}\rightarrow (\Aff\affi{A})^*$. We already know (lemma \ref{projetor_pim}) that this is a linear projector that transforms $v_P$ into $z_P$ and $v_m$ into $0$. Therefore:
	$$\alpha_0\cdot z_{R_0}+\ldots+\alpha_n\cdot z_{R_n}+\beta \cdot 0=0.$$
	
	But we already know (proposition \ref{base_zP}) that $z_{R_0},\ldots,z_{R_n}$ are a linear basis of $(\Aff\affi{A})^*$, hence concluding $\alpha_0=\ldots=\alpha_n=0$. Finally as $v_m\neq 0$ we must also conclude $\beta=0$.
	
	The elements $v_{R_0},\ldots,v_{R_n},v_m$ are then linearly independent. Using the dimension we conclude that they form a basis of $\Aff\Quad_m\affi{A}$.

	Regarding the statement on the dual basis, recall that affine functions are seen as elements on $\widehat{\Quad_m\affi{A}}$ using the immersion $\dd_m^*$ defined in (\ref{define_d_dt}). It suffices then to observe:
	$$\langle \delta^m_{\mathcal{R}},v_{R_j}\rangle=\delta^m_{\mathcal{R}}(R_j)=0,\quad \langle \delta^m_{\mathcal{R}},v_m\rangle=1,$$
	and using lemma \ref{projetor_pim}:	
	$$\langle \dd_m^* w_k, v_{R_j}\rangle=\langle w_k, \dd_m v_{R_j}\rangle=\langle w_k, z_{R_j}\rangle=\delta_{kj},$$ $$\langle \dd_m^*w_k,v_m\rangle=\langle w_k,\dd_m v_m\rangle=\langle w_k,0\rangle=0.$$
\end{proof}
\begin{prop}\label{propobasedelta}
	If $m$ is a non-degenerate metric, any affine referential $\mathcal{R}=(R_0,\ldots, R_n)$ on $\affi{A}$ determines a basis  $(\delta^m_{R_0},\ldots, \delta^m_{R_n},u_m)$ on $\widehat{\Quad_m\affi{A}}$.
\end{prop}
\begin{proof}
We know from proposition \ref{basevR} that $v_{R_0},\ldots, v_{R_n},v_m$ is a basis on the space $\Aff\Quad_m\affi{A}$, hence a system of linear coordinates on $\widehat{\Quad_m\affi{A}}$. To prove that the given elements form a basis it suffices to prove that they are linearly independent. Following $\langle v_{R_i},\delta^m_{R_j}\rangle =\delta^m_{R_j}(R_i)=\frac12 g(\vect{R_iR_j},\vect{R_iR_j})=\frac12 D_{ij}$, $\langle v_{R_i},u_m\rangle=1$,  $\langle v_m,\delta^m_{R_j}\rangle=1$, and $\langle v_m,u_m\rangle=0$, we only need to prove that the following matrix is non-degenerate:
$$\left[ \begin{matrix} \frac12D_{ij} & \one_n^t\\ \one_n & 0\end{matrix}\right].$$

However we already know from proposition \ref{nondegeneracy} that non-degeneracy of $g$ implies that the following matrix is non-degenerate:
$$\left[ \begin{matrix} D_{ij} & \one_n^t \\ \one_n & 0\end{matrix}\right].$$

Taking a product with a factor 2 in the first rows and a factor 1/2 in the last column, we complete our proof. 
\end{proof}
Observe that $\widehat{\Quad_m\affi{A}}$ and $\Aff\Quad_m\affi{A}$ are dual spaces of each other, but the basis $(\delta^m_{R_0},\ldots, \delta^m_{R_n},u_m)$ of the first space is not dual to the basis $(v_{R_0},\ldots, v_{R_n},v_m)$ of the other one. In particular, if we use $\delta_{R_i}^m,u_m$ as linear coordinate functions on $\Aff\Quad_m\affi{A}$, then $v_m$ has associated coordinate vector $[1\,\ldots\,1\,0]^t\in M_{(n+2)1}(\RR)$. If we use $v_{R_i},v_m$ as linear coordinate functions on $\widehat{\Quad_m\affi{A}}$, then the unit element $u_m$ has associated coordinate vector $[1\,\ldots\,1\,0]^t\in M_{(n+2)1}(\RR)$.

\begin{define}\label{CMcoordinates}
	For any fixed non-degenerate metric $m\in \Met\affi{A}$ and any affine referential $\mathcal{R}=(R_0,\ldots,R_n)$ on some $n$-dimensional affine space $\affi{A}$, we call the functions $\delta^m_{R_i}\in\Quad_m\affi{A}$ the Cayley-Menger o.d.system of quadratic coordinates (or CM coordinates) on $\affi{A}$.
\end{define}
\begin{remark}
Following proposition \ref{propobasedelta}, when $m\in\Met\affi{A}$ is non-degenerate the mapping $(\delta^m_{R_0},\ldots, \delta^m_{R_n},u_m)\colon\Aff\Quad_m\affi{A}\rightarrow \RR^{n+2}$ is a system of linear coordinates. It induces a system  $(\delta^m_{R_0},\ldots,\delta^m_{R_n}) \colon \Aff_1\Quad_m\affi{A}\rightarrow \RR^{n+1}$ of affine coordinates on the affine subspace of normalized functionals. Its composition with the quadratic immersion $v\colon \affi{A}\rightarrow \Aff_1\Quad_m\affi{A}\rightarrow \RR^{n+1}$ becomes an o.d.system of quadratic coordinates (definition \ref{odsystem}), called in the euclidean case half-squared distance coordinate functions. For any non-degenerate metric we call them Cayley-Menger coordinates.

These Cayley-Menger coordinates are quadratic functions and are overdetermined. Not any sequence of values is a valid sequence of CM coordinates for a point. Observe moreover that stating that $P$ has Cayley-Menger coordinate vector $(\delta_0,\ldots, \delta_n)$ in a given referential is the same as stating that $v_P$ has coordinate vector $[\delta_0\ldots\delta_n 1]^t$ in the system of linear coordinates $\delta^m_{R_0},\ldots, \delta^m_{R_n},u_m$ of $\Aff\Quad_m\affi{A}$. Moreover, using the linear projector $\dd_m\colon v_P\mapsto z_P$, one recovers the barycentric coordinate vector of $P$ as a linear combination of $(\delta_0,\ldots,\delta_n,1)$, hence as an affine expression of Cayley-Menger coordinate vector $(\delta_0,\ldots, \delta_n)$. The quadratic nature of the coordinates is then not a major problem, and ba\-ry\-centric coordinates can be recovered linearly from CM coordinates, without the need of square root computations. 

In modern technological applications where the position of robotic components is determined by relative distance measures obtained by interferometry or  sender-receiver delay of signal, the basic position information is that of relative distance, which can be better expressed in Cayley-Menger quadratic coordinates than in a non-existing fixed cartesian or affine referential \cite{cao2006sensor,gower1985properties,dokmanic2015euclidean,crippen1988distance,thomas2005revisiting}.
\end{remark}
Observe that from (\ref{d_on_functionals}):
$$\bar{v}\in \Aff_0\Quad_m\affi{A}\Rightarrow \langle\dd_m \bar{v},u\rangle=0\Rightarrow \dd_m \bar{v}\in \vect{A}\subset(\Aff\affi{A})^*.$$

As basic example we have $v_Q-v_P\in \Aff_0\Quad_m\affi{A}=\langle u_m\rangle^\circ$ (for points $P,Q\in\affi{A}$), because it factors by any addition of constants. As $\dd_m(v_P)=z_P$ (lemma \ref{projetor_pim}), we get  $\dd_m(v_Q-v_P)=z_Q-z_P=\vect{PQ}\in\vect{A}$. 
We have the exact sequence:
	\begin{equation} \label{umvmA}	\begin{tikzcd}
	0 \arrow{r}& \langle v_m\rangle \arrow[hookrightarrow]{r}{i}  & \Aff_0\Quad_m\affi{A} \arrow{r}{\dd_m} & \vect{A} \arrow{r}& 0.
	\end{tikzcd}
	\end{equation}
Observe for any $f\in\Aff\affi{A}$ and $\bar{v}_0\in\Aff_0\Quad_m\affi{A}=\langle u_m\rangle^\circ$, that $\dd_m\bar{v}_0$ is, by definition, obtained using any choice $\delta\in\Quad_m\affi{A}$ as:
$$\langle \dd f,\dd_m\bar{v}_0\rangle=\bar{v}_0(\delta+f)-\bar{v}_0(\delta).$$

On the other hand, by (\ref{define_d_dt}) we may consider $\Aff\affi{A}\subset\widehat{\Quad_m\affi{A}}$, and use $\bar{v}_0$ as a linear form on this space. We have
$$\bar{v}_0(\delta+f)-\bar{v}_0(\delta)=\langle \delta+f,\bar{v}_0\rangle-\langle \delta,\bar{v}_0\rangle=\langle f,\bar{v}_0\rangle.$$
Hence
\begin{equation}\label{double_d}\langle \dd f,\dd_mH\rangle=\langle f,H\rangle,\qquad \forall f\in\Aff\affi{A},\quad \forall H\in\Aff_0\Quad_m\affi{A}.\end{equation}

\begin{prop}
Consider a metric  $m\in\Met\affi{A}$ on the affine space $\affi{A}$, with hessian principal component $g\in S^2\vect{A}^*$. For any symmetric bilinear form $\CM$ on $\Aff\Quad_m\affi{A}$ the following conditions are equivalent:
	\begin{itemize}
		\item[1.] The restriction of $\CM$ to $\langle u_m\rangle^\circ=\Aff_0\Quad_m\affi{A}$ is projectable by $\dd_m$ as $-g$;
		\item[1'. ] For some affine referential $\mathcal{R}=(R_0,\ldots, R_n)$ of $\affi{A}$ there holds:
		$$\CM(v_{R_i}-v_{R_0},v_{R_j}-v_{R_0})=-g(\vect{R_0R_i},\vect{R_0R_j}),$$
		$$\CM(v_m,v_{R_i}-v_{R_0})=0,\quad \CM(v_m,v_m)=0; $$
		\item[1''. ] For any affine referential $\mathcal{R}=(R_0,\ldots, R_n)$ of $\affi{A}$ the previous conditions hold.
	\end{itemize}
\end{prop}
\begin{proof}
For the given referential $\mathcal{R}$, $v_{R_i}-v_{R_0},v_m$ form a a basis for $\langle u_m\rangle^\circ=\Aff_0\Quad_m\affi{A}$.
	
	For $\CM$ to restrict as a projectable symmetric bilinear form on $\langle u_m\rangle^\circ/\langle v_m\rangle$, the element $v_m$ should belong to the nullspace of $\CM$ restricted to $\langle u_m\rangle^\circ$, a condition that is encoded in the conditions $\CM(v_m,\cdot)=0$.
	
	The difference $v_{R_i}-v_{R_0}$ projects by $\dd_m$ as $z_{R_i}-z_{R_0}=\vect{R_0R_i}$. As we know, these vectors form a basis of $\vect{\affi{A}}$, which concludes our proof. 
\end{proof}
An interesting remark now is that, for any $P\in\affi{A}$ there holds $\langle v_P,u_m\rangle=1\neq 0$ and hence there exists a decomposition:
$$\Aff\Quad_m\affi{A}=\langle v_P\rangle\oplus \langle u_m\rangle^\circ.$$

We shall now consider a particular symmetric bilinear form on $\Aff\Quad_m\affi{A}$.

\begin{theorem}\label{teoremaprincipal}
For any non-null metric $m\in\Met\affi{A}$  with hessian principal component $g\in S^2\vect{A}^*$ on the affine space $\affi{A}$, there exists a unique symmetric bilinear form $\CM_m$ defined on $\Aff\Quad_m\affi{A}$ such that:
	\begin{enumerate}
		\item[1.] The restriction of $\CM_m$ to $\langle u_m\rangle^\circ=\Aff_0\Quad_m\affi{A}$ is projectable by $\dd_m$ to $-g$
		$$\CM_m(v_1,v_2)=-g(\dd_m(v_1),\dd_m(v_2)),\quad\forall v_1,v_2\in\Aff_0\Quad_m\affi{A}; $$
		\item[2.] For each point $P\in\affi{A}$ the functional $v_P\in \Aff_1\Quad_m\affi{A}$ is isotropic with respect to $\CM_m$:
		$$\CM_m(v_P,v_P)=0,\quad \forall P\in\affi{A}.$$
	\end{enumerate}
Moreover under the assumption (1.), the remaining condition (2.) is equivalent to any of the following:
	\begin{itemize}
	\item[2'a.] For each point $P\in\affi{A}$ there holds $\CM_m(v_P)=\delta^m_{P}$;
	\item[2'b.] At some point $P\in\affi{A}$ there holds $\CM_m(v_P)=\delta^m_{P}$;
	\item[2''a.] For each affine referential $\mathcal{R}=(R_0,\ldots, R_n)$ in $\affi{A}$ the elements $v_{R_i}$ are isotropic and $\CM_m(v_{R_i},v_m)=1$;
	\item[2''b.] For some affine referential $\mathcal{R}=(R_0,\ldots, R_n)\in\affi{A}$ the elements $v_{R_i}$ are isotropic and one of the points has $\CM_m(v_{R_0},v_m)=1$.
\end{itemize}
\end{theorem}
\begin{proof}
	We begin with conditions (1)+(2'b).
	
	$\bullet$ Existence and unicity: As $\Aff\Quad_m\affi{A}=\langle v_P\rangle\oplus\langle u_m\rangle^\circ$, any element of this space can be written in a unique way as $\alpha v_P+v$ (with $\alpha\in\RR$ and $v\in\langle u_m\rangle^\circ$)
	
	The given conditions:
	$$\CM_m\text{ symmetric bilinear},$$
	$$\CM_m(v_P)=\delta^m_P,$$
	$$\CM_m=-\dd_m^*(g)\text{ on }\langle u_m\rangle^\circ,$$
	univocally determine the bilinear form. Namely, the conditions imply that $\CM_m(v_P,v_P)=\delta^m_{P}(P)=0$ and it must be:
	$$\CM_m(\alpha v_P+v,\bar\alpha v_P+\bar v)=\alpha \langle \bar v,\delta^m_{P}\rangle + \bar\alpha\langle v,\delta^m_{P}\rangle - g(\dd_m(v),\dd_m(\bar v)) \quad v,\bar{v}\in \langle u_m\rangle^\circ,$$
	which is symmetric and satisfies all our conditions.

	For the equivalent characterizations we present next a proof that uses two cyclic arguments (2'b)$\Rightarrow$ (2)$\Rightarrow$(2'a)$\Rightarrow$(2'b), and (2''b)$\Rightarrow$(2)$\Rightarrow$(2''a)$\Rightarrow$(2''b), under the hypothesis (1).
	
	$\bullet$ Equivalence $(2'b)\Rightarrow (2)\Rightarrow (2'a)\Rightarrow (2'b)$ under assumption $(1)$

	Let us prove (1)+(2'b)$\Rightarrow$ (1)+(2)
	
	We know $\CM_m(v_P)=\delta^m_{P}$ for some specific point $P\in\affi{A}$. Hence:
	$$\begin{aligned}
		\CM_m&(v_Q,v_Q)=\CM_m(v_P+(v_Q-v_P),v_P+(v_Q-v_P))=\\
		&= \CM_m(v_P,v_P)+2\CM_m(v_P,v_Q-v_P)+\CM_m(v_Q-v_P,v_Q-v_P)=\\
		&=0+2\left(\delta^m_{P}(Q)-\delta^m_{P}(P)\right)-g(\vect{PQ},\vect{PQ})=\\
		&=0+g(\vect{PQ},\vect{PQ})+0-g(\vect{PQ},\vect{PQ})=0,
	\end{aligned}
	$$
	and we conclude that all elements $v_Q$ are isotropic for $\CM_m$.
	
	Let us prove (1)+(2)$\Rightarrow$ (1)+(2'a).
	
	If all elements $v_Q$ are isotropic, we have:
	$$-2\CM_m(v_P,v_Q)=\CM_m(v_P-v_Q,v_P-v_Q)=-g(\vect{PQ},\vect{PQ}),$$
	hence $\CM_m(v_P,v_Q)=\frac12g(\vect{PQ},\vect{PQ})=\delta^m_{P}(Q)$. This proves that $\CM_m(v_P)\circ v=\delta^m_P$, a m-quadratic function. Following lemma \ref{lemadeltav}, as $m\neq 0$, we conclude $\langle v_m,\CM_m(v_P)\rangle=1$, therefore $\CM_m(v_P)$ is in the hyperplane $\Quad_m\affi{A}\subset \widehat{\Quad_m\affi{A}}$ and represents precisely the function $\delta^m_P$.

	The implication (1)+(2'a)$\Rightarrow$ (1)+(2'b) is trivial. We have completed the equivalence of (2), (2'a), (2'b) under the assumption (1).
	
	$\bullet$ Equivalence $(2''b)\Rightarrow (2)\Rightarrow (2''a)\Rightarrow (2''b)$ under assumption $(1)$

	We prove now (1)+(2''b)$\Rightarrow$ (1)+(2'b) (which is equivalent to (1)+(2))

		Take the specific referential indicated in (1)+(2''b), for which:
		$$\CM_m(v_{R_0},v_m)=1,\quad \CM_m(v_{R_i},v_{R_i})=0.$$
		
		We have then:
		$$2\CM_m(v_{R_0},v_{R_j})=-\CM_m(v_{R_j}-v_{R_0},v_{R_j}-v_{R_0})=g(\vect{R_0R_j},\vect{R_0R_j})=2\delta^m_{R_0}(R_j).$$
		If we combine these two properties:
		$$\CM_m(v_{R_0},v_m)=1=\langle \delta^m_{R_0},v_m\rangle,$$
		$$\CM_m(v_{R_0},v_{R_j})=\langle\delta^m_{R_0},v_{R_j}\rangle,$$
		taking into account proposition \ref{basevR} we conclude that $\CM_m(v_{R_0})=\delta^m_{R_0}\in\Quad_m\affi{A}\subset \widehat{\Quad_m\affi{A}}$, so at least for some point we get the property given in (2'b), which, as was previously proven, together with (1) implies (2).

	Let us prove next (1)+(2)$\Rightarrow$ (1)+(2''a).
	
	One aspect is already given by (2), all $v_{R_i}$ are isotropic. The other aspect is to prove that $\CM_m(v_m,v_{R})=1$. By (2'a), we know that $\CM_m(v_{P})\circ v=\delta^m_P\circ v$, a $m$-quadratic function, therefore following lemma \ref{lemadeltav} $\CM_m(v_m,v_P)=\langle v_m,\CM_m(v_P)\rangle=1$, as we wanted to prove, for any point.
	
	Finally the implication (1)+(2''a)$\Rightarrow$ (1)+(2''b) is trivial.
\end{proof}	
\begin{define}\label{defineCM}
	The symmetric bilinear form $\CM_m$ on $\Aff\Quad_m\affi{A}$ characterized in theorem \ref{teoremaprincipal} shall be called Cayley-Menger product associated to the metric $m\in\Met\affi{A}$.
	
	Using the affine immersion $\Quad_m\affi{A}\hookrightarrow \Quad\affi{A}$ we have an induced linear projector $\widehat{\affi{A}}^2=(\Quad\affi{A})^*\rightarrow \Aff\Quad_m\affi{A}$ which allows to view Cayley-Menger product $\CM_m$ as a symmetric bilinear form on the quadratic hull $\widehat{\affi{A}}^2$, projectable by this mapping. This shall be called Cayley-Menger product associated to $m\in\Met\affi{A}$ on the quadratic hull $\widehat{\affi{A}}^2$.
\end{define}
\begin{remark}
	As stated in remark \ref{metricanula}, in the case of the null metric $m\in\Met\affi{A}$ we must not mistake $\delta^m_P$ (null quadratic function) in $\widehat{\Quad_m\affi{A}}$ with the null element of this vector space. Conditions [1.] and [2'a.] in theorem \ref{teoremaprincipal} do determine a unique (non-null) Cayley-Menger bilinear form $\CM_0$ also in this case, a bilinear form that satisfies all our conditions. However condition [2.] is not equivalent to the other ones, as for example the null bilinear form on $\Aff\Quad_m\affi{A}$ would also satisfy [1.] and [2.] without ever taking any value $\delta^m_P\in\Aff_1\Quad_m\affi{A}$. As $m=0$ is a rather uninteresting case for applications we may always assume $m\neq 0$. The ``void'' Cayley-Menger (degenerate) bilinear form $\CM_0$ however might be an interesting tool in some situations, in the absence of any metric structure on $\affi{A}$.
\end{remark}
\begin{prop}\label{propCMref}
	The following properties hold, for the Cayley-Menger product $\CM_m$ associated to a metric $m\in\Met\affi{A}$:
		\begin{enumerate}
		\item On the subspace $\Aff_0\Quad_m\affi{A}$ there holds $\dd\circ\CM_m=-g\circ\dd_m$.
		\item For any referential $\mathcal{R}=(R_0,\ldots, R_n)$ of $\affi{A}$, $\CM_m$ transforms the basis $(v_{R_0},\ldots, v_{R_n},v_m)$ of $\Aff\Quad_m\affi{A}$ into the elements $(\delta^m_{R_0},\ldots, \delta^m_{R_n},u_m)$ of $\widehat{\Quad_m\affi{A}}$ (which by proposition \ref{propobasedelta} is a basis of this vector space if $m$ is non-degenerate).
		\item For any pair of affine referentials $\mathcal{R}=(R_0,\ldots, R_n)$ and $\mathcal{S}=(S_0,\ldots, S_n)$ on $\affi{A}$ the matrix representation of the linear morphism $\CM_m\colon \widehat{\Quad_m\affi{A}}\rightarrow \Aff\Quad_m\affi{A}$ using $(v_{R_0},\ldots, v_{R_n},v_m)$ as basis on the first space and linear coordinates $(v_{S_0},\ldots, v_{S_n},v_m)$ on the second space is the following (mixed reference Cayley-Menger matrix):
		 $$\left[\begin{matrix} \frac12D^{\mathcal{R}\mathcal{S}}&\one^t \\ \one & 0\end{matrix}\right]\qquad D^{\mathcal{R}\mathcal{S}}_{ij}= g(\vect{R_iS_j},\vect{R_iS_j}).$$
		\item The Gram matrix associated to $\CM_m$ using $(v_{R_0},\ldots, v_{R_n},v_m)$ as basis is the following (Cayley-Menger matrix):
		$$\left[\begin{matrix} \frac12D^{\mathcal{R}}&\one^t \\ \one & 0\end{matrix}\right]\qquad D^{\mathcal{R}}_{ij}= g(\vect{R_iR_j},\vect{R_iR_j}).$$
		\item If $P,Q$ have barycentric coordinate vectors $p,q \in\WW_n$ in the referential $\mathcal{R}$ then their squared pseudodistance is determined by:
		$$d^2(P,Q)=-\frac12(q-p)^t\cdot D^{\mathcal{R}}\cdot (q-p).$$
		\end{enumerate}
\end{prop}
\begin{proof}
For any functional $\bar{v}_0\in\Aff_0\Quad_m\affi{A}$ we want to prove that the covectors $-\dd\circ \CM_m \bar{v}_0$ and $g\circ\dd_m\bar{v}_0$ coincide. It suffices to apply duality with any vector $\vect{x}\in\vect{\affi{A}}$. As we know from (\ref{umvmA}) that all these vectors have the form $\vect{x}=\dd_m\bar{w}$ for some choice of $\bar{w}\in\Aff_0\Quad_m\affi{A}$. In this situation:
$$\langle -\dd\circ\CM_m \bar{v}_0,\dd_m \bar{w}\rangle=\langle -\CM_m \bar{v}_0,\bar{w}\rangle = \langle g\circ \dd_m \bar{v}_0,\dd_m \bar{w}\rangle,$$
where we use formula (\ref{double_d}) and the fact that $\CM_m$ on $\Aff_0\Quad_m\affi{A}$ projects by $\dd_m$ as $-g$.

Hence both covectors are coincident on elements of the form $\dd_m\bar{w}$ and we conclude $g\circ\dd_m=-\dd\circ\CM_m$.

Property (2'a) in Theorem \ref{teoremaprincipal} shows that $\CM(v_{R_i})=\delta^m_{R_i}$. To prove that $\CM_m(v_m)$ equals $u_m$ it suffices to observe from  the definition of $u_m$ that $\langle u_m,v_{R_i}\rangle=1$ and $\langle u_m,v_m\rangle=0$, while $\CM_m(v_m,v_{R_i})=1$ and $\CM(v_m,v_m)=0$, following (2''a) and (1') in the same theorem.

When we use two referentials $\mathcal{R}$, $\mathcal{S}$ as given in the statement we get:
$$\langle \CM_m(v_{R_i}),v_{S_j}\rangle =\langle \delta^m_{R_i},v_{S_j}\rangle=\frac12\delta^m_{R_i}(S_j)=\frac12 g(\vect{R_iS_j},\vect{R_iS_j}),$$
$$\CM_m(v_{R_i},v_m)=\langle \delta^m_{R_i},v_m\rangle=1,\quad \CM_m(v_m,v_m)=\langle u_m,v_m\rangle=0.$$

The matrix associated to $\CM_m$ using the given pair of referentials is then the one given in the statement, and the Gram matrix on the next part of our statement corresponds to the case of a single referential $\mathcal{S}=\mathcal{R}$ on the quadratic hull.

Finally, if we call $g$ the symmetric bilinear form on $\vect{A}$ associated to the metric $m$ we observe from property (1.) in the definition of $\CM_m$:
$$d^2(P,Q)=\langle z_P-z_{Q},g(z_P-z_Q)\rangle=-\langle v_P-v_{Q},\CM_m(v_P-v_{Q})\rangle,$$
because $v_P-v_{Q}$ belongs to $\langle u_m\rangle^\circ$.
Using a coordinate representation with respect to the basis induced by $\mathcal{R}$:
$$d^2(P,Q)=-\left[ \alpha_0\ldots \alpha_n\beta\right]\cdot \left[\begin{matrix} \frac12D^{\mathcal{R}}&\one^t \\ \one & 0\end{matrix}\right]\cdot \left[ \alpha_0\ldots \alpha_n\beta\right]^t.$$
Where $\alpha_0,\ldots, \alpha_n,\beta$ represent the coordinates of $v_P-v_{Q}$ in the basis $v_{R_0},$ $\ldots,$ $v_{R_n},$ $v_m$. As $v_Q$ projects to $\widehat{\affi{A}}$ as $z_Q$, the components $\alpha_0,\ldots,\alpha_n$ are the linear coordinates of $z_P-z_{Q}$, hence they are the difference $p-q$ of the corresponding barycentric coordinate vectors. Moreover, as $\one_n\cdot p=\one_n\cdot q=1$, we conclude $\one_n\cdot \alpha=0$, and the expression above becomes independent of the component $\beta$:
$$d^2(P,Q)=-\frac12(q-p)^t\cdot D^{\mathcal{R}}\cdot (q-p).$$
This formula is in consonance with the already know (\ref{qpDelta}) taking into account that $\delta^m_P(Q)=\frac12 d^2(P,Q)$ and the already stated relations $-\frac12D^\mathcal{R}=G=2\Delta$ for our matrix representations.
\end{proof}

\begin{remark}
	The existence of a bilinear form on certain spaces, with properties analogous to those described by theorem \ref{teoremaprincipal} is known in the literature, sometimes with an ad-hoc construction and sometimes with a more intrinsic description. Theorem 5.5 in \cite{bertram2004linear} is a good example, for the case of a non-degenerate metric $m\in\Met\affi{A}$. We may translate this theorem in our language, saying that a bilinear form $G$ is determined on a certain vector space (generated by quadratic functions $\delta^m_P$) using as a property that $G(\delta^m_P,\delta^m_Q)=\frac12\|\vect{PQ}\|^2$. As we shall see in the following, this is in fact the inverse quadratic form associated to our Cayley-Metric quadratic form, for the case of non-degenerate metrics.
	
	Also in certain applications \cite{crippen1988distance} authors consider a space constructed as the free affine span of a set of points (called amalgamation space associated to this set), where distances determine a bilinear form (called by the authors as Schoenberg's quadratic form \cite{schoenberg1935remarks}). However this presentation is not functorial and heavily relies on a choice of points. 
\end{remark}
We may use Cayley-Menger bilinear form to relate two affine retractions of the quadratic mappings $v\colon \affi{A}\rightarrow \Aff\Quad_m\affi{A}$ and $\delta^m\colon \affi{A}\rightarrow \widehat{\Quad_m\affi{A}}$:
\begin{lemma}\label{immersproj}
	Consider the mappings $v\colon \affi{A}\rightarrow \Aff_1\Quad_m\affi{A}$ and $\delta^m\colon \affi{A}\rightarrow \Quad_m\affi{A}$ determined by a choice of metric $m\in\Met\affi{A}$. For any pair of points $P,Q\in\affi{A}$ and 1-dimensional weight $(\alpha,\beta)\in\WW_1$ there holds:
	\begin{itemize}
		\item $v(\alpha P+\beta Q)-\alpha v(P)-\beta v(Q)\in\langle v_m\rangle=\RR$,
		\item $\delta^m(\alpha P+\beta Q)-\alpha\delta^m(P)-\beta\delta^m(Q)\in\langle u_m\rangle=\RR$.		
	\end{itemize}
	The constant value, for both cases is $\frac{-\alpha\beta}{2}\cdot\CM_m(v_P,v_Q)$.
\end{lemma}
\begin{proof}
	We consider first $\bar v=v(\alpha P+\beta Q)-\alpha v(P)-\beta v(Q)$.

	For any $\delta\in\Quad_m\affi{A}$ and any $f\in\Aff\affi{A}$ there holds:
	$$\langle \bar v,\delta+f\rangle=\langle \bar v,\delta\rangle+\langle v_{\alpha P+\beta Q},f\rangle-\alpha\langle v_{P},f\rangle-\beta\langle v_{Q},f\rangle.$$
	
	For affine $f$, we know $$\langle v_{\alpha P+\beta Q},f\rangle=f(\alpha P+\beta Q)=\alpha f(P)+\beta f(Q)=\alpha\langle v_P,f\rangle+\beta\langle v_Q,f\rangle.$$
	Therefore $\langle \bar v,\delta+f\rangle=\langle \bar v,\delta\rangle$ for any affine function $f\in\Aff\affi{A}$. Recall that the director vector space associated to $\Quad_m\affi{A}$ is $\Aff\affi{A}$. Therefore $\bar v$ is an affine operator that takes constant value on $\Quad_m\affi{A}$. As $\Quad_m\affi{A}$ is an hyperplane on $\widehat{\Quad_m\affi{A}}$, determined by $v_m(\delta)=1$, this implies that $\bar v\in\langle v_m\rangle$, which is the first part of our statement.
	
	This implies now that $\CM_m(\bar v)\in\langle \CM_m(v_m)\rangle$, which is our second statement: $$\delta^m(\alpha P+\beta Q)-\alpha\delta^m(P)-\beta\delta^m(Q)\in \langle u_m\rangle.$$
	To determine the constant value in both statements, observe that following corollary \ref{coroCMimmersion}, $\delta^m$ is a quadratic mapping, therefore:
	$$\begin{aligned} \delta^m&(\alpha P+\beta Q)-\alpha\delta^m(P)-\beta\delta^m(Q)=\alpha(\alpha-\beta)\delta^m(P)+\beta(\beta-\alpha)\delta^m(Q)+\\&+4\alpha\beta\delta^m\left(\frac{P+Q}{2}\right)-\alpha(\alpha+\beta)\delta^m(P)-\beta(\alpha+\beta)\delta^m(P)=\\
		&=-2\alpha\beta\left( \delta^m(P)+\delta^m(Q)-2\delta^m\left(\frac{P+Q}{2}\right)\right) \end{aligned}$$
	We know that this is a constant function (depending on the choice of $P,Q,\alpha,\beta$). Taking value at the point $(P+Q)/2$, for example, and knowing from definition of $\delta^m$ that $\delta^m_P((P+Q)/2)=\delta^m_Q((P+Q)/2)=\frac14\delta^m_P(Q)$ we get:
	$$\begin{aligned} \delta^m(\alpha P+\beta Q)&-\alpha\delta^m(P)-\beta\delta^m(Q)=-2\alpha\beta\cdot \frac14\delta^m_P(Q)=-\frac12 \alpha\beta\delta^m_P(Q)=\\
		&=\frac{-\alpha\beta}{2}\CM_m(v_P,v_Q)\end{aligned}$$
\end{proof}
As a consequence composition of the quadratic mapping $\delta^m$ with the affine projector $p\colon \Quad_m\affi{A}\rightarrow \Quad_m\affi{A}/\RR$ determines an affine mapping $p\circ\delta^m\colon \affi{A}\rightarrow \Quad_m\affi{A}/\RR$.
We may also say that $v\colon \affi{A}\rightarrow \Aff_1\Quad_m\affi{A}$ composed with the projector $\dd_m$ is an affine mapping. However in this case, this composition is clearly the mapping $z\colon \affi{A}\hookrightarrow\widehat{\affi{A}}$, which we already knew to be affine.

\begin{prop}\label{propdeltasQ}
Consider a metric $m\in\Met\affi{A}$ with associated hessian principal component $g\in S^2\vect{A}^*$.	For any fixed point $Q\in\affi{A}$ consider the gradient at $Q$ of quadratic functions, an affine mapping $\nabla_Q\colon\Quad_m\affi{A}\rightarrow \vect{A}^*$ with a natural linear extension to the linear hull $\widehat{\Quad_m\affi{A}}$. Consider the induced dual mapping $\nabla_Q^*\colon \vect{A}\rightarrow \Aff\Quad_m\affi{A}$.  These mappings determine a commutative diagram:
	\begin{equation*}
	\begin{tikzcd}
		\vect{A} \arrow{rr}{\nabla_Q^*} \arrow{d}{-g}  & & \Aff\Quad_m\affi{A} \arrow{d}{\CM_m} \\
		\vect{A}^* & & \widehat{\Quad_m\affi{A}} \arrow{ll}{\nabla_Q}.
	\end{tikzcd}
\end{equation*}
Moreover the nullspace of $\nabla_Q$ is generated by $u_m$ and $\delta^m_Q$, and  the image of $\nabla_Q^*$ is the set of affine functionals that factor by $\RR$ and vanish at $\delta^m_Q\in\Quad_m\affi{A}$.
\end{prop}
\begin{proof}
Observe that $\nabla_Q$ factors by $\RR$ (addition of a constant does not change the gradient of a quadratic function, at any point). Therefore $\nabla_Q^*$ takes values on the subspace $\Aff_0\Quad_m\affi{A}$.

Recall that the gradient at any point of any affine function is the linear principal component of the affine function. Consequently for any director vector $\vect{a}\in\vect{A}$ and any affine function $f\in\Aff\affi{A}$ and following (\ref{double_d}) there holds:
$$\langle \dd_m(\nabla_Q^*(\vect{a})),\dd f\rangle=\langle f,\nabla_Q^*(\vect{a})\rangle=\langle \vect{a},\nabla_Q f\rangle =\langle \vect{a},\dd f\rangle,$$
hence $\dd_m\circ \nabla_Q^*(\vect{a})=\vect{a}$ for any director vector $\vect{a}\in\vect{A}$. 

Recall that by definition the restriction of $\CM_m$ to $\Aff_0\Quad_m\affi{A}$ is the pull-back of $-g$ by $\dd_m\colon \Aff_0\Quad_m\affi{A}\rightarrow \vect{A}$. We have then a commutative diagram:
	\begin{equation*}
	\begin{tikzcd}
		\vect{A} \arrow{rr}{\nabla_Q^*}  & & \Aff_0\Quad_m\affi{A} \arrow{d}{\CM_m} \arrow{rr}{\dd_m}&& \vect{A} \arrow{d}{-g}\\
		\vect{A}^* & & \widehat{\Quad_m\affi{A}} \arrow{ll}{\nabla_Q} && \vect{A}^*, \arrow{ll}{\dd_m^*}
	\end{tikzcd}
\end{equation*}
which, considering that $\dd_m\circ \nabla_Q^*=\Id$ implies that $-g\colon \vect{A}\rightarrow \vect{A}^*$ closes the diagram on the left hand side.

For the second part of the statement, on  $\Aff\affi{A}\subset\widehat{\Quad_m\affi{A}}$ the gradient mapping is simply the computation of the principal component. Therefore $\nabla_Q$ is surjective onto $\vect{A}^*$ and by dimension its nullspace is a 2-dimensional subspace.

As $\nabla_Q$ factors by $\RR$, the element $u_m\in\widehat{\Quad_m\affi{A}}$ is in the nullspace of $\nabla_Q$. Moreover, by definition $\delta^m_Q$ has null gradient at $Q$, hence $\delta^m_Q$ is also in this nullspace. We know $\langle \delta^m_Q,v_m\rangle=1$ hence $\delta^m_Q$ is linearly independent with $u_m$. Therefore the nullspace of $\nabla_Q$ is generated by these two elements, and the image of $\nabla_Q^*$ is the set of affine operators that vanish on $u_m$ (hence belong to $\Aff_0\Quad_m\affi{A})$ and at the same time vanish on $\delta^m_Q$.
\end{proof}

\begin{prop}\label{RRAdual}
Consider a metric $m\in\Met\affi{A}$ with hessian principal component $g\in S^2\vect{A}^*$ and Cayley-Menger bilinear form $\CM_m\colon \Aff\Quad_m\affi{A}\rightarrow \widehat{\Quad_m\affi{A}}$. Consider the bilinear form $\text{swap}(a,b)=(b,a)$ defined for $(a,b)\in\RR\oplus\RR$.
The following is a commutative diagram, where the horizontal arrows are isomorphisms adjoint to each other:

	\begin{equation*}
	\begin{tikzcd}
		\RR\oplus\RR\oplus\vect{A}  \arrow{d}{(\text{swap},-g)} \arrow{rr}{v_m\oplus v_Q\oplus \nabla_Q^*}  & & \Aff\Quad_m\affi{A} \arrow{d}{\CM_m} \\
		\RR\oplus\RR\oplus\vect{A}^* & & \widehat{\Quad_m\affi{A}}. \arrow{ll}{(v_m,v_Q,\nabla_Q)}
	\end{tikzcd}
	\end{equation*}

Hence $\widehat{\Quad_m\affi{A}}$ with the Cayley-Menger metric is isometric to the direct sum of the standard hyperbolic plane $(\RR^2,\text{swap})$ and the metric vector space $(\vect{A},-g)$.
\end{prop}
\begin{proof}
By linearity, we only need to study these linear mappings on elements 	$(1,0,\vect{0})$, $(0,1,\vect{0})$ and $(0,0,\vect{a})$ of $\RR\oplus\RR\oplus\vect{A}$.

Firstly, as $u_m$ has vanishing hessian principal component and has vanishing gradient and value 1 at any point $Q$:	\begin{equation*}
	\begin{tikzcd}
	(1,0,\vect{0}) \arrow{rr}{v_m\oplus v_Q\oplus \nabla_Q^*}&& v_m \arrow{r}{\CM_m}& u_m \arrow{rr}{(v_m,v_Q,\nabla_Q)}&& (0,1,0).
	\end{tikzcd}
\end{equation*}
Secondly, as $\delta^m_Q$ has $g$ as hessian principal component, and  takes vanishing value and gradient at $Q$:
\begin{equation*}
	\begin{tikzcd}
		(0,1,\vect{0}) \arrow{rr}{v_m\oplus v_Q\oplus \nabla_Q^*}&& v_Q \arrow{r}{\CM_m}& \delta^m_Q \arrow{rr}{(v_m,v_Q,\nabla_Q)}&& (1,0,0).
	\end{tikzcd}
\end{equation*}
Finally, from proposition \ref{propdeltasQ}, we know that $\CM_m\nabla_Q^*(\vect{a})$ has gradient $-g(\vect{a})$ at point $Q$. Moreover, as $\CM_m$ is self-adjoint, we observe that $$\langle v_m,\CM_m\nabla_Q^*(\vect{a})\rangle=\langle u_m,\nabla_Q^*(\vect{a})\rangle=0$$
and $$\langle v_Q,\CM_m\nabla_Q^*(\vect{a})\rangle=\langle \delta^m_Q,\nabla_Q^*(\vect{a})\rangle=\langle \nabla_Q\delta^m_Q,\vect{a}\rangle=0,$$ and we conclude:
\begin{equation*}
	\begin{tikzcd}
		(0,0,\vect{a}) \arrow{rr}{v_m\oplus v_Q\oplus \nabla_Q^*}&& \nabla_Q^*(\vect{a}) \arrow{r}{\CM_m}& \CM_m\nabla_Q^*(\vect{a}) \arrow{rr}{(v_m,v_Q,\nabla_Q)}&& (0,0,-g(\vect{a})).
	\end{tikzcd}
\end{equation*}
Thus proving our result.

Regarding the surjectivity of $(v_m,v_Q,\nabla_Q)$, we may observe that it transforms $\delta^m_Q$ into $(1,0,\vect{0})$, and also transforms any affine function $f\in\Aff\affi{A}\subset \widehat{\Quad_m\affi{A}}$ into $(0,f(Q),\dd f)$. By dimension computation $(v_m,v_Q,\nabla_Q)$ and $v_m\oplus v_Q\oplus \nabla_Q^*$ will be isomorphisms. Using the natural duality pairing of $\RR\oplus\RR\oplus\vect{A}^*$ with $\RR\oplus\RR\oplus\vect{A}$, clearly $(v_m,v_Q,\nabla_Q)$ is adjoint to $v_m\oplus v_Q\oplus \nabla_Q^*$, thus completing the proof.
\end{proof}
\begin{remark}
	Following proposition \ref{RRAdual}, elements of $\widehat{\Quad_m\affi{A}}$ can be seen in a simpler fashion as elements in $\RR\times\RR\times \vect{A}^*$, when one fixes an origin $O$  on the affine space (the referential point $Q$ of this proposition).	
	This mapping represents a function $\delta^m_P-\frac12r^2\in\Quad_m\affi{A}$ as $\left(1,\frac12\|\vect{OP}\|^2-\frac12r^2,-g(\vect{OP})\right)\in\RR\times\RR\times\vect{A}^*$. It also represents an affine function $f(X)=\omega(\vect{OX})+c\in\Aff\affi{A}$ as $(0,c,\omega)=(0,f(O),\dd f)$.
	
	In Möbius geometry (see \cite{cecil2008lie} for example) on an $n$-dimensional space $\affi{A}$ with euclidean structure $g\in S^2\vect{A}^*$ each hypersphere is represented as a line on a certain $n+2$-dimensional Lorentz metric space. The immersion is constructed using a particular point and stereographic projection. The representation given by Möbius takes a hypersphere with center $P$ and radius $r$ (which we may identify with a function $\delta^m_P-\frac12r^2$) into a line generated by  $$(u,v,\vec{x})=\left(\frac{1+\|\vect{OP}\|^2-r^2}{2},\frac{1-\|\vect{OP}\|^2+r^2}{2},\vect{OP}\right)\in\RR^2\oplus\vect{\affi{A}}.$$
	In the space $\RR^2_{u,v}\oplus \vect{\affi{A}}$ Möbius geometry considers a Lorentzian structure using $(\dd u^2-\dd v^2)\oplus g$. It suffices to consider the mapping $(u,v,\vec{x})\mapsto (u+v,\frac{u-v}{2},-g(\vec{x}))\in\RR^2\oplus \vect{A}^*$ together with the isomorphisms given in proposition \ref{RRAdual}
	to  see that Möbius ad-hoc representation of hyperspheres is just a particular representation of $\widehat{\Quad_m\affi{A}}$ with Cayley-Menger metric (which is a Lorentz metric, if $g$ is euclidean). In this representation hyperspheres get identified with $m$-quadratic functions $\delta^m_P-\frac12r^2\in\Quad_m\affi{A}\subset \widehat{\Quad^m\affi{A}}$.
	
This immersion into a Lorentz space and its projectivization was also used by Pedoe \cite{pedoe1937representation} in his study of a product between circles of the plane. In our formalism, this Pedoe product is just the application of Cayley-Menger metric to some specific elements of $\widehat{\Quad_m\affi{A}}$. 
	
	We get a formulation of the classical Möbius (and also Lie Sphere) geometry in an intrinsic fashion, in the light of the natural Cayley-Menger bilinear form presented in theorem \ref{teoremaprincipal}. The implications of our results for polarities on projective spaces will be our focus in a companion paper.
\end{remark}
\begin{corollary}
For any point $Q\in\affi{A}$, the vector subspace $\vect{H}_Q=\langle v_Q,v_m\rangle\subset \Aff\Quad_m\affi{A}$ with the restricted Cayley-Menger bilinear form $\CM_m$ is a hyperbolic plane (signature $(1,1,0)$), and its $\CM_m$-orthogonal complement $\vect{H}_Q^\bot$ is isometric to $(\vect{A},-g)$ by $\nabla_Q^*\colon\vect{A}\rightarrow \vect{H}_Q^\bot\subseteq \Aff\Quad_m\affi{A}$. 
\end{corollary}

\begin{corollary}\label{signaturaCM}
	If $m$ is a metric on $\affi{A}$ with  signature $(\pi,\nu,\rho)$, then $\CM_m$ is a symmetric bilinear form on $\Aff\Quad_m\affi{A}$ with signature $(\nu+1,\mu+1,\rho)$
\end{corollary}
\begin{corollary}\label{nondegenerateCM}
	If $m$ is a non-degenerate metric on $\affi{A}$ then $\CM_m$ is a non-degenerate metric on $\Aff\Quad_m\affi{A}$.
\end{corollary}

We are now in the situation to characterize points of the Cayley-Menger quadric (the image of the quadratic immersion $A\subset \Aff\Quad_m\affi{A}$) as the intersection of the cone of $\CM_m$-isotropic vectors with an affine hyperplane orthogonal to $v_m$:
\begin{theorem}\label{CMquadricAff}
Consider for any non-null metric $m\in\Met\affi{A}$ the Cayley-Menger bilinear form $\CM_m$   on $\Aff\Quad_m\affi{A}$. For any element $\bar v\in\Aff\Quad_m\affi{A}$ and for the quadratic immersion $v\colon P\in\affi{A}\mapsto v_P\in\Aff\Quad_m\affi{A}$ there holds: $$\bar{v}\in\Img v\Leftrightarrow \CM_m(\bar v,\bar v)=0,\quad \CM_m(\bar v,v_m)=1.$$
\end{theorem}
\begin{proof}
	One implication is immediate. If $\bar v=v_P$ we know:
	$$\CM_m(\bar v,\bar v)=\CM_m(v_P,v_P)=\langle v_P,\delta^m_{P}\rangle=\delta^m_{P}(P)=0,$$
	$$\CM_m(\bar v,v_m)=\CM_m(v_P,v_m)=\langle v_m,\delta^m_{P}\rangle=1.$$

	Le us prove the converse. Take $z=\dd_m(\bar v)\in\widehat{\affi{A}}=(\Aff{\affi{A}})^*$ for the natural projector $\dd_m\colon\Aff\Quad_m\affi{A}\rightarrow \widehat{\affi{A}}$ studied in lemma \ref{projetor_pim}. Observe that the nullspace of this projection is $v_m$.

	We assume now $\langle u_m,\bar v\rangle=\CM_m(v_m,\bar v)=1$. Hence recalling that $u_m=u\circ\dd_m$ we have $\langle u,z\rangle=1$. We know that this implies $z=z_P\in\affi{A}\subset\widehat{\affi{A}}$ for some point $P\in\affi{A}$.
	
	Let us prove that $\bar v=v_P$ using $\CM_m(\bar v,\bar v)=0$.
	
	Being $z=z_P\in\widehat{\affi{A}}$ and as $v_m$ generates the nullspace of $\dd_m$ we may conclude $\bar v=v_P+c\cdot v_m\in\Aff\Quad_m\affi{A}$ for some element $c\in\RR$.
	
	$$\begin{aligned} 0&=\CM_m(\bar v,\bar v)=\CM_m(v_P+c\cdot v_m,v_P+c\cdot v_m)=\\
		&=\CM_m(v_P,v_P)+2c\cdot \CM_m(v_m,v_P)+c^2\cdot \CM_m(v_m,v_m)=2c,\end{aligned}$$
	hence $c=0$ and we conclude $\bar v=v_P$.
\end{proof}
\begin{corollary}\label{invCMquadric}
	For a non-degenerate metric, the associated Cayley-Menger bilinear form $\CM_m\colon \Aff\Quad_m\affi{A}\rightarrow \widehat{\Quad_m\affi{A}}$ is invertible, and the Cayley-Menger quadratic mapping $\delta^m\colon\affi{A}\rightarrow \Quad_m\affi{A}$ is an immersion (corollary \ref{coroCMimmersion}). The image of $\delta^m$ is characterized by:
	$$\delta\in\Img \delta^m\Leftrightarrow\CM_m^{-1}(\delta,\delta)=0,\quad \CM_m^{-1}(\delta,u_m)=1.$$
\end{corollary}
\begin{proof}
Following propositions \ref{propCMref} and \ref{propobasedelta}, the image of $\CM_m$ is the whole space $\widehat{\Quad_m\affi{A}}$. Moreover we know $v\circ\CM_m=\delta^m$ (property 2. in theorem \ref{teoremaprincipal}) and how $\CM_m$ characterizes the image of $v$ in theorem \ref{CMquadricAff}. Hence $\delta$ has the form $\delta^m_P$ if and only if $\CM_m^{-1}(\delta)$ has the form $\CM_m^{-1}(\delta^m_P)=v_P$. Therefore:
$$\delta\in\Img \delta^m\Leftrightarrow \CM_m\left(\CM_m^{-1}(\delta),\CM_m^{-1}(\delta)\right)=0,\quad \CM_m\left(\CM_m^{-1}(\delta),v_m\right)=0.$$

Taking into account that $u_m=\CM_m(v_m)$ this can be written as:
$$\CM_m^{-1}(\delta,\delta)=0,\quad \CM_m^{-1}(\delta,u_m)=1.$$
\end{proof}
\begin{remark}
Consider a metric $m\in\Met\affi{A}$ with Gram matrix representation $G=G^m_{\mathcal{R}}$ in the referential $\mathcal{R}=(R_0,\ldots, R_n)$. We know its relation to the squared distance matrix $D_{ij}=\widehat{g}(\vect{R_iR_j},\vect{R_iR_j})$ by $G=-\frac12 D$.	Following proposition \ref{propCMref} the bilinear form $\CM_m$ has Gram matrix in the basis $v_{R_0},\ldots, v_{R_n},v_m$ given by 
\begin{equation} \label{matrizCayleyMenger}\left[\begin{matrix} -G&\one^t\\ \one & 0\end{matrix}\right].
\end{equation}
When $m$ is non-degenerate $\CM_m$ has an inverse and if we use $\delta^m_{R_i}=\CM_m(v_{R_i})$ and $u_m=\CM_m(v_m)$ as basis on $\widehat{\Quad_m\affi{A}}$ then Gram representation of $\CM_m^{-1}$ on this basis is again (\ref{matrizCayleyMenger}), but if we use as basis the dual basis of $v_{R_i},v_m$, then Gram representation of $\CM_m^{-1}$ is the inverse of (\ref{matrizCayleyMenger}).
\end{remark}

The presentation of Cayley-Menger bilinear form in definition \ref{defineCM} was characterized following intrinsic properties enumerated in theorem \ref{teoremaprincipal}. All these properties use spaces and objects that are functorially derived from the metric affine space $(\affi{A},m)$. One would expect a covariant behaviour of Cayley-Menger bilinear form $\CM_m$ with respect to affine transformations.

Consider any affine transformation $\varphi\colon \affi{B}\rightarrow \affi{A}$. Composition of affine mappings with quadratic mappings is quadratic. We have then induced mappings:
	\begin{itemize}
		\item[-] Linear mapping $\varphi_{\Quad{}}\colon \Quad\affi{A}\rightarrow \Quad\affi{B}$ defined by $(\varphi_{\Quad{}} \delta) (P)=\delta(\varphi(P))$ for any $P\in\affi{B}$ and $\delta\in\Quad\affi{A}$.
		\item[-] Linear mapping $\varphi_{\Aff{}}\colon \Aff{\affi{A}}\rightarrow \Aff{\affi{B}}$, restriction of the previous one to the subspace $\Aff{\affi{A}}\subset\Quad\affi{A}$ and taking values in $\Aff\affi{B}$ (composition of affine functions is an affine function).
		\item[-] Linear mapping $\varphi_{\Met{}}\colon \Met\affi{A}\rightarrow \Met\affi{B}$, induced by $\varphi_{\Quad{}}$ from the quotient space $\Met\affi{A}=\Quad\affi{A}/\Aff\affi{A}$ to $\Met\affi{B}=\Quad\affi{B}/\Aff\affi{B}$.
		\item[-] Affine mapping $\varphi_{\Quad{}}\colon \Quad_m\affi{A}\rightarrow \Quad_{\bar{m}}\affi{B}$, restriction of $\varphi_{\Quad{}}$ to the subspace $\Quad_m\affi{A}$, that takes values in $\Quad_{\bar{m}}\affi{B}$ with $\bar{m}=\varphi_{\Met{}}m$,  due to the definition of $\varphi_{\Met{}}m$.
		\item[-] Linear mapping $\varphi_{\Quad{}}^*\colon \Aff\Quad_{\varphi_{\bar{m}}}\affi{B}\rightarrow \Aff\Quad_m\affi{A}$, obtained by composition with $\varphi_{\Quad{}}$.
		\item[-] Linear mapping $\widehat{\varphi}\colon \widehat{\affi{B}}\rightarrow\widehat{\affi{A}}$, dual to $\varphi_{\Aff{}}$.
		\item[-] Linear mapping $\dd\varphi\colon \vect{B}\rightarrow \vect{A}$, restriction of the previous one to $\vect{B}$.
	\end{itemize}
\begin{lemma}\label{lemaphivP} For any affine mapping $\varphi\colon\affi{B}\rightarrow \affi{A}$, for any point $P\in\affi{B}$ and for the unit functions $u_A\in\Quad\affi{A}$, $u_B\in\Quad\affi{B}$ there holds::
$$\varphi_{\Quad{}}^* v_P=v_{\varphi(P)},\quad \varphi_{\Quad{}}u_A=u_B$$
\end{lemma}
\begin{proof}
	For any $\delta\in\Quad\affi{A}$ we have:
	$$\langle \varphi^*_{\Quad{}}v_P,\delta\rangle=\langle v_P,\varphi_{\Quad{}}(\delta)\rangle=(\varphi_{\Quad{}}\delta)(P)=\delta(\varphi(P))$$
	hence $\varphi^*_{\Quad{}}v_P$ takes on any quadratic function $\delta$ precisely the value of $\delta$ at the point $\varphi(P)$. This is precisely the definition of $v_{\varphi(P)}$.
	
	Regarding the unit function it is evident that $u_A(\varphi(P))=1$ for any $P\in\affi{B}$, hence $\varphi_{\Quad{}}^*u_A=u_B$.
\end{proof}
As a consequence of $\varphi_{\Quad{}}u_A=u_B$ we deduce that $\varphi_{\Quad{}}^*$ transforms $\RR$-projectable affine operators on $\Quad_{\bar{m}}\affi{A}$ into $\RR$-projectable affine operators on $\Quad_{m}\affi{B}$ and as on these spaces $\dd^A_m$ and $\dd^B_{\bar{m}}$ given in (\ref{define_d_dt}) are simply the computation of the principal linear components, we get a commutative diagram:
	\begin{equation}\label{comuta}
	\begin{tikzcd}
		\Aff_0\Quad_m\affi{A}   \arrow{rr}{\dd_m^A}  & & \vect{A}  \\
		\Aff_0\Quad_{\bar{m}}\affi{B} \arrow{u}{\varphi_{\Quad{}}^*}  \arrow{rr}{\dd_{\bar{m}}^B}& & \vect{B}. \arrow{u}{\dd\varphi}
	\end{tikzcd}
	\end{equation}
\begin{lemma}
	If $g_m$ is the Hessian principal component of the metric $m\in\Met\affi{A}$, then the Hessian principal component of the metric $\bar{m}=\varphi_{\Met{}}m$ is $\dd\varphi^*g_m$.
\end{lemma}
\begin{proof}
Take any quadratic function $\delta\in\Quad\affi{A}$ representing $m\in\Met\affi{A}$. The metric $\varphi_{\Met{}}m$ is represented by the quadratic function $\varphi_{\Quad{}}\delta$

The quadratic function $\delta$ is associated to a bilinear quadratic form $\widehat{g}$ on $\widehat{\affi{A}}$ given in terms of $\delta$ by (\ref{formulagtildezetas}). Moreover for any $P,Q\in\affi{B}$ we know that $z_{\varphi(P)}=\widehat{\varphi}(z_P)$. Hence:
$$\widehat{g}(\widehat{\varphi}(z_P),\widehat{\varphi}(z_Q))=\widehat{g}(z_{\varphi(P)},z_{\varphi(Q)})=4\delta\left(\frac{\varphi(P)+\varphi(Q)}{2}\right)-\delta(\varphi(P))-\delta(\varphi(Q))$$
In the same way $\varphi_{\Quad{}}\delta$ is a quadratic function represented by another bilinear form:
$$\begin{aligned}\widehat{g}_2(z_P,z_Q)&=4(\varphi_{\Quad{}}\delta)\left(\frac{P+Q}{2}\right)-(\varphi_{\Quad{}}\delta)(P)-(\varphi_{\Quad{}}\delta)(Q)=\\
	&=4\delta\left(\varphi\left(\frac{P+Q}{2}\right)\right)-\delta\left(\varphi(P)\right)-\delta\left(\varphi(Q)\right)\end{aligned}$$
Taking into account that $\varphi$ is affine, we conclude that $\widehat{g}_2=\widehat{\varphi}^*\widehat{g}$. As both Hessian principal components are obtained by restriction of these bilinear forms to $\vect{A}\subseteq\widehat{A}$ and $\vect{B}\subseteq\widehat{B}$, and as the restriction of $\widehat{\varphi}\colon\widehat{\affi{B}}\rightarrow \widehat{\affi{A}}$ is the linear mapping $\dd\varphi\colon\vect{B}\rightarrow \vect{A}$ we complete our proof.
\end{proof}

\begin{theorem}\label{functorialCM}
	For any affine mapping $\varphi\colon \affi{B}\rightarrow \affi{A}$ from some affine space $\affi{B}$ to a metric affine space $(\affi{A},m)$ (where $m\in\Met\affi{A})$, such that the induced metric $\bar{m}=\varphi_{\Met{}}m\in\Met\affi{B}$ is non-null, the Cayley-Menger bilinear form $\CM_m$ associated to $m$ and the Cayley-Menger bilinear form $\CM_{\bar{m}}$ associated to $\bar{m}$ satisfy the relation:
	$$\CM_{\bar{m}}(\bar{v},\bar{w})=\CM_m(\varphi_{\Quad{}}^*\bar{v},\varphi_{\Quad{}}^*\bar{w}),\,\forall \bar{v},\bar{w}\in\Aff\Quad_{\bar{m}}\affi{B}$$
\end{theorem}
\begin{proof}
	It suffices to prove that the bilinear form defined on elements $\bar{v},\bar{w}\in\Aff\Quad_{\bar{m}}\affi{B}$ by $\CM_m(\varphi_{\Quad{}}^*\bar{v},\varphi_{\Quad{}}^*\bar{w})$ fulfills conditions (1) and (2) imposed for $\CM_{\bar{m}}$ in theorem \ref{teoremaprincipal}. We have to prove that:
	$$\CM_m(\varphi_{\Quad{}}^*\bar{v}_1,\varphi_{\Quad{}}^*\bar{v}_2)=-g_{\bar{m}}(\dd^B_{\bar{m}} \bar{v}_1, \dd^B_{\bar{m}} \bar{v}_2),\quad \forall \bar{v}_1,\bar{v}_2\in\Aff_0\Quad_{\bar{m}}\affi{B}$$
	$$\CM_m(\varphi_{\Quad{}}^*v_P,\varphi_{\Quad{}}^*v_P)=0,\qquad \forall P\in\affi{B}$$
	
	The second condition is immediate from lemma \ref{lemaphivP} because $\CM_m$ satisfies condition (2) of theorem \ref{teoremaprincipal}:
	$$\CM_m(\varphi_{\Quad{}}^*v_P,\varphi_{\Quad{}}^*v_P)=\CM_m(v_{\varphi(P)},v_{\varphi(P)})=0$$
	
	For the first condition we observe that the hessian principal component associated to $\bar{m}$ is $\dd\varphi^* g_m$ (where $g_m$ is Hessian principal component associated to $m\in\Met\affi{A}$, hence:
	$$-g_{\bar{m}}(\dd^B_{\bar{m}}\bar{v}_1,\dd^B_{\bar{m}}\bar{v}_2)= -g_{m}(\dd\varphi\circ\dd^B_{\bar{m}}\bar{v}_1,\dd\varphi\circ\dd^B_{\bar{m}}\bar{v}_2)=\ldots$$
	using commutative diagram (\ref{comuta}) we get:
	$$\ldots=-g_{m}(\dd^A_{m}\circ\varphi_{\Quad{}}^*\bar{v}_1,\dd^A_{m}\circ\varphi_{\Quad{}}^*\bar{v}_2)=\ldots$$
	and using now property (1) from theorem \ref{teoremaprincipal} that defined $\CM_m$:
	$$\ldots=\CM_m(\varphi_{\Quad{}}^*\bar{v}_1,\varphi_{\Quad{}}^*\bar{v}_2)$$
	which completes our proof
\end{proof}

\begin{corollary}
	If $\varphi\colon\affi{A}\rightarrow \affi{A}$ is an affine isometry on the metric affine space $(\affi{A},m)$, then $\varphi_{\Quad{}}^*$ is a linear isometry of the vector space $\Aff\Quad_m\affi{A}$, with respect to Cayley-Menger bilinear form $\CM_m$. If $\CM_m$ is invertible, then the dual mapping $\varphi_{\Quad{}}$ is an affine transformation on $\Quad_m\affi{A}$ extending to a linear transformation on the linear hull $\widehat{\Quad_m\affi{A}}$, which represents an isometry of $\widehat{\Quad_m\affi{A}}$ with respect to Cayley-Menger inverse bilinear form $\CM_m^{-1}$.
\end{corollary}
%
%
%
%

\begin{remark}\label{diagrama}
We may illustrate the mappings relating all objects described in this section, in the case $\dim \affi{A}=1$ (affine line) with a diagram given in figure \ref{figura}, where double arrows represent quadratic mappings between affine spaces and single arrows represent affine mappings (that are determined by linear mappings on the corresponding linear hulls). On $\Aff{\Quad_m\affi{A}}$ we are representing the isotropic cone of Cayley-Menger product $\CM_m$, and on the dual space $\widehat{\Quad_m\affi{A}}$ the isotropic cone associated to the inverse Cayley-Menger product $\CM_m^{-1}$ (which exists if $m$ is non-degenerate). 
\begin{figure}\caption{Cayley-Menger bilinear form and natural morphisms.}\label{figura}
	\centering\includegraphics[angle=90,origin=c,width=0.59\linewidth]{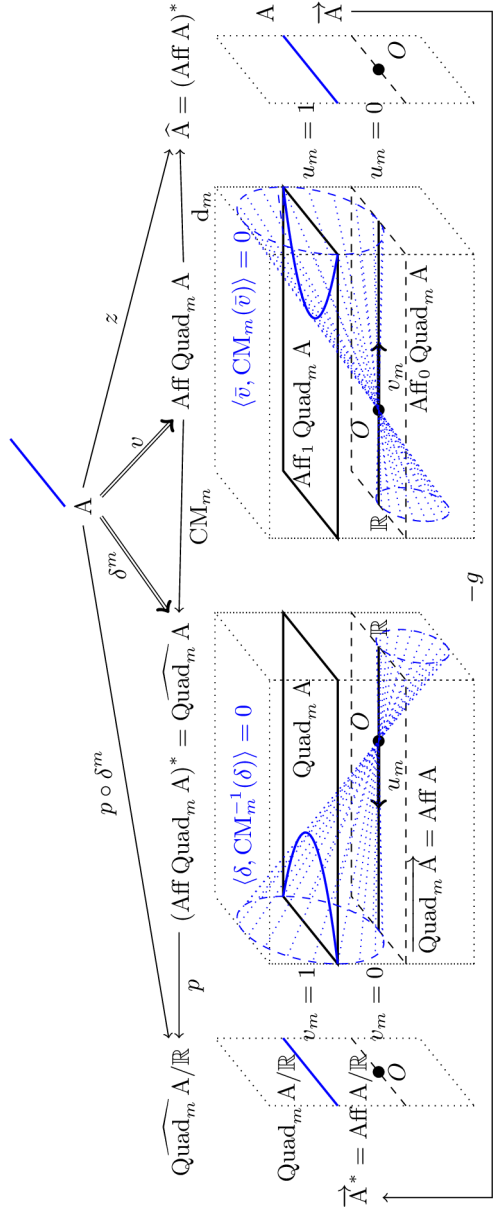}\end{figure}
	Observe that, even though $v\colon \affi{A}\hookrightarrow \Aff\Quad_m\affi{A}$ is a quadratic immersion, it has a linear retraction $\dd_m\colon \Aff\Quad\affi{A}\rightarrow \widehat{\affi{A}}$. This linear retraction, maps the affine subspace $\left(u_m=1\right)\equiv \Aff_1\Quad_m\affi{A}\subset\Aff\Quad_m\affi{A}$ into the affine subspace $\left(u=1\right)\equiv\affi{A}\subset \widehat{\affi{A}}$. In the same manner, for a non-degenerate metric, even though $\delta^m$ is a quadratic mapping, $p$ acts as a linear retraction if we identify the affine space $\Quad_m\affi{A}/\RR$ (of $m$-covector fields on $\affi{A}$)  with the affine space $\affi{A}$,  mapping each $m$-covector field to the unique point where the associated gradient covector vanishes.
\end{remark}

\section{Example}\label{secejemplo}
We will illustrate the applicability of all the tools presented in the previous sections in a simple situation. For a larger variety of applications of Cayley-Menger matrix (hence of our specific presentation as bilinear form) we suggest the consideration of \cite{havel1991some}.

Consider a 2-dimensional real affine space $\affi{A}$ and a given referential $\mathcal{R}=(R_0,R_1,R_2)$. We may represent points on $\affi{A}$ using its barycentic coordinate vector (a weight), and director vectors of $\vect{A}$ using its hollow coordinate vector. For example, for the midpoint $R_{02}=(R_0+R_2)/2\in\affi{A}$, for $P=\inv_{R_0}R_1=2R_0-R_1\in\affi{A}$ inversion of $R_1$ with respect to $R_0$ (see definition \ref{simetrico}), and for the director vector $\vect{R_0R_2}\in\vect{A}$ we have the matrix representations $r_{02},p\in\WW_2$, $x_{02}\in\HH_2$:
$$R_{02}\to r_{02}=[\,1/2\,0\,1/2\,]^t,\quad P\to p=[\,2\,-1\,0\,]^t,\quad \vect{R_0R_2}\to x_{02}=[-1\,0\,1\,]^t$$
Consider a specific quadratic function $\delta\in\Quad\affi{A}$ characterized by a matrix $S$ of its values $s_{ij}=\delta\left(\frac{R_i+R_j}{2}\right)$ at the referential and its midpoints, as indicated below. Using (\ref{matrizDelta}) we obtain its barycentric coordinate representation $\Delta$:
$$S=\left[ \begin{matrix} 1&4&-2 \\ 4&9&-1 \\ -2&-1&1\end{matrix}\right]\Rightarrow \Delta=2S-\frac12\left((\diag S)^t\one+\one^1(\diag S)\right)=\left[\begin{matrix} 1&3&-5\\3&9&-7 \\-5&-7&1\end{matrix}\right]$$
The value $\delta(P)$ can be obtained using barycentric coordinates of the point as $\delta(P)=p^t\cdot \Delta\cdot p$. The affine function $f\in\Aff\affi{A}$ that has at $R_0,R_1,R_2$ the same values as $\delta$ is characterized by the row coefficient vector $c=\diag S=\diag\Delta=[\,1\,9\,1\,]$. The homogeneous component of $\delta$ at $R_{02}$, at $P$ and the $\mathcal{R}$-reduced component of $\delta$ are new quadratic functions $\delta^h_{R_{02}}$, $\delta^h_P$, $\delta^0_{\mathcal{R}}$ with barycentric representations given by (\ref{formula4b}) and (\ref{coefdelta0}):
\begin{equation}\label{variasdelta}\Delta^h_{R_{02}}=\left[ \begin{matrix} 3&5&-3\\5&11&-5\\-3&-5&3\end{matrix}\right],\quad\Delta^h_P=\left[ \begin{matrix} 4&8&0\\8&16&0\\0&0&8\end{matrix}\right],\quad \Delta^0_{\mathcal{R}}=\left[ \begin{matrix} 0&-2&-6\\-2&0&-12\\-6&-12&0\end{matrix}\right]\end{equation}
Each of these functions $\delta^h_{R_{02}},\delta^h_P,\delta^0_{\mathcal{R}},\delta$ is associated, by theorem \ref{largepropo}, to a bilinear form on $\widehat{\affi{A}}$. In the basis $z_{R_0},z_{R_1},z_{R_2}$ of this vector space these bilinear forms have Gram matrix $G=2\Delta$, $G^h_{R_{02}}=2\Delta^h_{R_{02}}$, $G^h_P=2\Delta^h_P$, $G^0_{\mathcal{R}}=2\Delta^0_{\mathcal{R}}$, respectively. Simple computations show that these symmetric matrices don't share the same inertia indexes. For the first two ones (homogeneous at a point) this inertia is $(2,0,1)$. Following remark \ref{remark43} this implies that the hessian principal component has index $(2,0,0)$ and $\delta$ is convex.

Quadratic functions $\delta,\delta^h_{R_{02}},\delta^h_P,\delta^0_{\mathcal{R}}$ differ from each other by an affine function. They represent the same metric $m\in\Quad\affi{A}/\Aff\affi{A}$. The associated quadratic forms on $\widehat{\affi{A}}$ have the same restriction on the subspace $\vect{A}\subset \widehat{\affi{A}}$. Using $\vect{R_0R_1},\vect{R_0R_2}$ as a basis of $\vect{A}$, their restriction has Gram matrix:
$$\left[\begin{matrix} -1&-1\\1&0\\0&1\end{matrix}\right]^t\cdot G\cdot \left[\begin{matrix} -1&-1\\1&0\\0&1\end{matrix}\right]=\left[\begin{matrix} 8&-8\\-8&24\end{matrix}\right]\quad \text{(same for }G^h_{R_{02}},G^h_P,G^0_{\mathcal{R}}\text{)}$$
which is positive-definite. Bilinear forms on $\widehat{\affi{A}}$ associated to other quadratic functions representing the same metric, however, need not share the same inertia index.

 We have a unique representative $\delta^0_{\mathcal{R}}$ of the metric, that vanishes at the referential, and $G=2 \Delta^0_{\mathcal{R}}$ is the linear hull Gram matrix representation of this metric when we use $v_{R_0}$, $v_{R_1}$, $v_{R_2}$ as basis. 

Following proposition \ref{propobasedelta}, quadratic functions $\delta^m_{R_0},\delta^m_{R_1},\delta^m_{R_2},u_m$ determine a basis on $\widehat{\Quad_m\affi{A}}$. In this space we have the inverse Cayley-Menger metric $\CM_m^{-1}$, whose Gram matrix in this basis is given by (\ref{matrizCayleyMenger}):
\begin{equation}\label{CMexample}\CM_m\rightarrow\left[\begin{matrix} -G^0_{\mathcal{R}} & \one^t \\ \one^t & 0 \end{matrix}\right]=
\left[ \begin{matrix} 0&4&12&1\\4&0&24&1\\12&24&0&1 \\ 1&1&1&0\end{matrix}\right]\end{equation}
We would like to know if $\delta$ belongs to Cayley-Menger quadric, if it has the form $\delta^m_Q$ for some point.

The $m$-quadratic function $\delta$ takes values $1,9,1$ at points $R_0,R_1,R_2$, hence its coordinates in the dual basis of $v_{R_0},v_{R_1},v_{R_2},v_m$ is $d=[\,1\,9\,1\,1\,]$ and:
$$\CM_m^{-1}(\delta,\delta)=\left[\,1\,9\,1\,1\right]\cdot \left[ \begin{matrix} 0&4&12&1\\4&0&24&1\\12&24&0&1 \\ 1&1&1&0\end{matrix}\right]^{-1}\cdot \left[\begin{matrix} 1\\9\\1\\1\end{matrix} \right]=-4\neq 0$$
hence $\delta$ doesn't belong to Cayley-Menger quadric in this case, it can not be written as $\delta^m_Q$ for any point.

We know $\CM_m^{-1}(\delta,\delta)=-4$ and $\CM_m^{-1}(u_m,u_m)=0$. For any $c\in\RR$ we have $$\begin{aligned}\CM_m^{-1}&(\delta+cu_m,\delta+cu_m)=-4+c^2\cdot 0+2c\CM_m^{-1}(\delta,u_m)=\\&=-4+0+2c\cdot\langle \delta,v_m\rangle=-4+2c\end{aligned}$$
We conclude that $\delta+2u_m$ belongs to the Cayley-Menger quadric, it has the form $\delta^m_Q$. There exists a point $Q\in\affi{A}$ and a constant $r=2$ such that $\delta=\delta^m_Q-\frac12r^2$. The function $\delta$ we chose for our example has the form $\delta^m_Q-\frac12r^2$ and it vanishes on an $m$-sphere centered at some point $Q$ with radius $r=2$. We may deduce the position of this point. As $\delta^m_Q=\delta+2$ takes values $1+2$, $9+2$, $1+2$ at the points $R_0,R_1,R_2$, respectively, using $v_{R_0},v_{R_1},v_{R_2},v_m$ as linear coordinates, the element $\delta^m_Q\in\Quad_m\affi{A}\subset\widehat{\Quad_m\affi{A}}$ has coordinates $(3,11,3,1)$. We know $v_Q=\CM_m^{-1}(\delta^m_Q)$ hence multiplication with the inverse of matrix (\ref{CMexample}) determines the coordinates of $v_Q$ in the basis $v_{R_0},v_{R_1},v_{R_2},v_m$:
$$v_Q=\frac12 v_{R_0}+0 v_{R_1}+\frac12 v_{R_2}-3v_m$$
projection by $\dd_m$ shows that $z_Q= 1/2 v_{R_0}+0 v_{R_1}+1/2 v_{R_2}$ and that the center of the $m$-sphere has barycentric coordinates $(1/2,0,1/2)$ in our referential, hence this center is precisely $R_{02}\in\affi{A}$. Indeed, we may consult our matrix expressions (\ref{variasdelta}) and find that $$\Delta^h_{R_{0,2}}-2\one^t\cdot \one=\Delta,$$ reflecting the fact that $\delta=\delta^h_{R_{02}}-2u$.

In arbitrary (pseudo)metric spaces one may use Cayley-Menger matrix with respect to an arbitrary referential to solve several metric questions using linear algebra techniques. Most applications \cite{alfakih2018euclidean,cao2006sensor,crippen1988distance,dokmanic2015euclidean,thomas2005revisiting,thomas2020clifford,liberti2014euclidean,gower1985properties,havel1991some} refer to the euclidean case but, as we see in this paper, using appropriate geometric concepts one extends this theory to arbitrary affine spaces or even to affine bundles, in the presence of a metric with any signature.

\section*{Appendix: Characterization of quadratic functions}
\renewcommand{\thesection}{Ap} In this appendix we shall prove certain results that show that the quadraticity condition (\ref{quadraticity}) characterizes functions on affine spaces whose affine coordinate representation follows the classical quadratic polynomial expressions. In order to understand these results as ``non-trivial'', recall that definition \ref{definequadratico} relies on barycentric combinations of points (\ref{quadraticity}) and not on any affine coordinate  polynomial representation.

\begin{lemma}
	If $\delta$, $\tilde \delta$ are quadratic  mappings between affine spaces $\affi{A}$ and $\affi{B}$ and coincident at two points $P,Q\in\affi{A}$ and at the corresponding midpoint $\frac{P+Q}{2}\in\affi{A}$, then they are also coincident at all the points of the affine line $\langle P,Q\rangle\subseteq\affi{A}$
\end{lemma}
\begin{proof}
	Points in the affine line have the form $\alpha\cdot P+\beta\cdot Q$ with $\alpha+\beta=1$. This lemma is then immediate consequence of quadraticity condition (\ref{quadraticity}).
\end{proof}

\begin{lemma}	
	(Homothety substitution rule) If $\delta$, $\tilde \delta$ are quadratic  mappings between affine spaces $\affi{A}$ and $\affi{B}$ and coincident at three points $P_0,P_1,P_2$ and at their midpoints $P_{ij}=\frac{P_i+P_j}{2}$, then for any real value $\alpha$ the quadratic functions $\delta_1,\delta_2$ are also coincident at the points $\bar{P_0}=P_0$, $\bar P_1=P_0+\alpha(P_1-P_0)$, $\bar P_2=P_2$ and at their midpoints $\bar P_{ij}=\frac{\bar P_i+\bar P_j}{2}$.
\end{lemma}
\begin{proof}
	We shall consider certain auxiliary points as shown in the following diagram:
	$$R_0=-\frac12 P_0+P_1+\frac12 P_2,\quad  R_1=\frac34 P_1+\frac14 P_2,\quad R_2=\frac32 P_0-\frac12 P_2$$
	\begin{center}
		\tikzstyle{Pontos_dados}=[fill=black, draw=black, shape=circle]
		\tikzstyle{Pontos_auxiliares}=[fill={rgb,255: red,191; green,191; blue,191}, draw=black, shape=circle]
		\tikzstyle{pontos_objetivo}=[fill=white, draw=black, shape=rectangle]
		\tikzstyle{Retas_dadas}=[-, fill=black]
		\tikzstyle{Retas_auxiliares}=[-, draw={rgb,255: red,128; green,128; blue,128}]
		\begin{tikzpicture}[scale=0.3]
				\node [style={Pontos_dados},label=$P_{11}$] (P11) at (0, 0) {};
				\node [style={Pontos_dados},label=$P_{01}$] (P01) at (6, 0) {};
				\node [style={Pontos_dados},label=$P_{00}$] (P00) at (12, 0) {};
				\node [style={Pontos_dados},label=$P_{12}$] (P12) at (6, 6) {};
				\node [style={Pontos_dados},label=$P_{02}$] (P02) at (12, 6) {};
				\node [style={Pontos_dados},label=$P_{22}$] (P22) at (12, 12) {};
				\node [style={pontos_objetivo},label=$\bar P_{01}$] (PP01) at (21, 0) {};
				\node [style={pontos_objetivo},label=$\bar P_{11}$] (PP11) at (30, 0) {};
				\node [style={pontos_objetivo},label=$\bar P_{02}$] (PP02) at (21, 6) {};
				\node [style={Pontos_auxiliares},label=$R_{0}$] (R0) at (0, 6) {};
				\node [style={Pontos_auxiliares},label=$R_{1}$] (R1) at (3, 3) {};
				\node [style={Pontos_auxiliares},label=$R_{2}$] (R2) at (12, -6) {};
				\draw [style={Retas_dadas}] (P11) to (P22);
				\draw [style={Retas_dadas}] (P00) to (P11);
				\draw [style={Retas_dadas}] (P00) to (P22);
				\draw [style={Retas_auxiliares}] (R0) to (R2);
				\draw [style={Retas_auxiliares}] (P00) to (R2);
				\draw [style={Retas_auxiliares}] (R0) to (PP02);
				\draw [style={Retas_auxiliares}] (P22) to (PP11);
				\draw [style={Retas_auxiliares}] (P00) to (PP11);
		\end{tikzpicture}
	\end{center}
	As $\delta,\tilde\delta$ are quadratic and coincident at $P_{00},P_{11}$ and at its midpoint $P_{01}$, we conclude that they are coincident on the affine line $\langle P_{00},P_{11}\rangle$, hence:
	$$\delta(\bar P_{11})=\tilde\delta(\bar P_{11}),\qquad \delta(\bar P_{01})=\tilde\delta(\bar P_{01})$$
	We still need to prove that they are coincident at $\bar{P}_{02}$.
	
	As $\delta,\tilde\delta$ are quadratic and coincident at $P_{00},P_{22}$ and at its midpoint $P_{02}$, we conclude that 
	$$\delta(R_{2})=\tilde \delta(R_{2})$$
	As $\delta,\tilde\delta$ are quadratic and coincident at $P_{11},P_{22}$ and at its midpoint $P_{12}$, we conclude that
	$$\delta(R_1)=\tilde\delta(R_1)$$
	
	Taking into account that $P_{01}$ is the midpoint of $R_0$ and $R_{2}$ and that $R_1=\frac34 R_0+\frac14 R_2$  is on the same line, the quadraticity condition (\ref{quadraticity}) for $\delta$ shows that:
	$$\begin{aligned}
		\delta(R_1)&=\frac34\left(\frac34-\frac14\right)\delta(R_0)+\frac14\left(\frac14-\frac34\right)\delta(R_{2})+4\frac34\frac14\delta(P_{01})\Rightarrow\\
		&\Rightarrow  
		\delta(R_0)=\frac{1}{3}\delta(R_{2})-2\delta(P_{01})+\frac83\delta(R_{1})\end{aligned}
	$$
	Also for $\tilde\delta$ the same argument and formula holds. Taking into account that $\delta,\tilde\delta$ are coincident at $R_1,R_2,P_{01}$ we conclude that $$\delta(R_0)=\tilde\delta(R_0)$$
	
	Finally, as $\delta,\tilde\delta$ are quadratic and coincident at $R_0$, $P_{02}$ and at its midpoint $P_{12}$, they are coincident on the line $\langle P_{12},P_{02}\rangle$, hence $$\delta(\bar P_{02})=\tilde \delta(\bar P_{02})$$
	which completes our proof.
\end{proof}

\begin{lemma}
	(General substitution rule) If $\delta,\tilde\delta$ are quadratic  mappings between affine spaces $\affi{A},\affi{B}$, coincident at points $P_0,P_1,\ldots,P_k$ and at the corresponding midpoints $P_{ij}=\frac12(P_i+P_j)$, then for any real number $x\in\RR$, at the points $\bar P_0=P_0,\quad \bar P_1=P_1,\quad \bar P_2=P_1+x(P_2-P_0)$, $\bar P_i=P_i$ ($\forall i> 2$) and at the corresponding midpoints both quadratic  mappings are also coincident.
\end{lemma}
\begin{proof}
	Consider the points $P_0,P_1,P_2,P_i$ ($i> 2$).
	
	Using the homothety substitution rule, with factor $2$, with the points $P_0,P_1$ and each of the remaining points of the original list, we may replace and take $\bar P_1=P_0+2(P_1-P_0)$. Both $\delta,\tilde\delta$ are then coincident at the points:
	$$P_0,P_0+2(P_1-P_0),P_2,P_i.\quad i>2$$
	and at the corresponding midpoints.
	
	Using the homothety substitution rule, with factor $2x$, with the points $P_0,P_2$ and each of the remaining points of the new list, we conclude that both $\delta,\tilde\delta$ are then coincident at the points:
	$$P_0,P_0+2(P_1-P_0),P_0+2x(P_2-P_0),P_i$$
	and at the corresponding midpoints.
	
	Using the homothety substitution rule with $\alpha=\frac12$ for the second and third point of this list, and with each of the remaining points, we conclude that both $\delta,\tilde\delta$ are coincident at the points:
	$$P_0,P_0+2(P_1-P_0),P_1+x(P_2-P_0),P_i$$
	and at the corresponding midpoints.
	
	Using finally the homothety substitution rule, with factor $1/2$ with the first two points, we conclude that $\delta,\tilde\delta$ are coincident at the points
	$$P_0,P_1,P_1+x(P_2-P_0),P_i$$
	and at the corresponding midpoints.
	
\end{proof}
The substitution rule above is the affine description of the elementary transformations common in linear algebra. From this property we finally conclude our main result:

\begin{prop} If $\delta,\tilde\delta$ are quadratic  mappings between affine spaces $\affi{A},\affi{B}$ and coincident at all points $P_0,P_1,\ldots,P_k$ and at the corresponding midpoints $P_{ij}=\frac{P_i+P_j}{2}$, then they are coincident at any affine combination $\alpha_0P_0+\alpha_1P_1+\ldots+\alpha_kP_k$ of these points (where $\alpha_0+\alpha_1+\ldots+\alpha_k=1$)
\end{prop}
\begin{proof}
	We begin with an application of the Homothety substitution.
	
	Being $\delta,\tilde\delta$ coincident at $P_0,P_1,\ldots, P_k$ and at the corresponding midpoints, the homothety substitution rule applied to $P_0,P_1,P_i$ with scalar $\alpha_1$ shows that they are also coincident at:
	$$P_0,\,\bar P_1=(1-\alpha_1)P_0+\alpha_1P_1,\,P_2,\ldots, P_k$$
	and at the corresponding midpoints.
	
	Our general substitution rule applied to the first three points and scalar $\alpha_2$ shows then that $\delta,\tilde\delta$ are also coincident at:
	$$P_0,\bar P_1,\, \bar P_2=(1-\alpha_1-\alpha_2)P_0+ \alpha_1P_1+\alpha_2P_2,\,P_3,\ldots, P_k$$
	and at the corresponding midpoints.
	
	Another application of the general substitution rule on the first, third and fourth points with scalar $\alpha_3$ shows that $\delta,\tilde\delta$ are also coincident at:
	$$P_0,\bar P_1,\bar P_2,\,\bar P_3=(1-\alpha_1-\alpha_2-\alpha_3)P_0+\alpha_1P_1+\alpha_2P_2+\alpha_3P_3,\,P_4,\ldots,P_k$$
	and at the corresponding midpoints.
	
	Iterating the same argument we finally obtain a point where $\delta$ and $\tilde\delta$ are coincident:
	$$(1-\alpha_1-\alpha_2-\ldots-\alpha_k)P_0+\alpha_1P_1+\ldots+\alpha_kP_k$$
	This point is the affine combination in our statement (recall $\alpha_0+\alpha_1+\ldots+\alpha_k=1$, hence the first term of this addition is $\alpha_0P_0$)
\end{proof}

\begin{theorem}
	Consider two affine spaces $\affi{A}$ and $\affi{B}$. Fix an affine referential $\mathcal{R}=(R_0,R_1,\ldots,R_n)$ and its midpoints $R_{ij}=\frac{R_i+R_j}{2}$. For any given $\binom{n+2}{2}$ points $Q_{ij}\in\affi{B}$ (where $Q_{ij}=Q_{ji}$, $0\leq i,j\leq n$), there exists a unique quadratic  mapping $\delta\colon \affi{A}\rightarrow \affi{B}$ such that $\delta(R_{ij})=Q_{ij}$. This quadratic  mapping can be given as:
	$$\delta(x_0 R_0+\ldots+x_nR_n)=\sum_{i,j} x_ix_j \Delta_{ij},\qquad \Delta_{ij}=2Q_{ij}-\frac12 Q_{ii}-\frac12 Q_{jj}\in\affi{B}$$
\end{theorem}
\begin{proof}
	The unicity relies on the previous lemmas. Any pair of quadratic  mappings that take common known values $Q_{ij}$ at points $R_{ij}$ must be the same on every affine combination of the points $R_i$.
	
	For the existence, it suffices to see that the proposed formula is a quadratic mapping with the given values at the mentioned points. 	Let us represent by $x$ the barycentric coordinate vector for any point $P$. There holds $\one\cdot x=1$, for the row matrix $\one=[1\ldots 1]$. 
	
	We may write the formula proposed in the statement as $$\delta(P)=x^t\cdot \Delta\cdot x$$ where $\Delta$ is a square matrix taking points $\Delta_{ij}\in\affi{B}$ as values.
	
	Observe that $\sum_{ij} x_ix_j=1$ when $\sum x_i=1$, hence the expression given for $\delta$ is an affine combination of points $\Delta_{ij}\in\affi{B}$. Observe also that for these points there holds $\Delta_{ij}=\Delta_{ji}$ and $\Delta_{ii}=Q_{ii}$. The proposed formula is a well-defined affine combination and clearly takes the given values $Q_{ij}$ at all points $R_{ij}=\frac12R_i+\frac12R_j$.

	We want to prove that our definition of $\delta$ is quadratic, that is, when we consider two points $R,S$ and consider $\alpha,\beta$ with $\alpha+\beta=1$ :
	\begin{equation}\label{quadracon}\delta(\alpha R+\beta S)=\alpha(\alpha-\beta)\delta(R)+ \beta(\beta-\alpha)\delta(S)+4\alpha\beta\delta\left(\frac{R+S}{2}\right)\end{equation}
	
	Take the barycentric coordinate vectors $x,y$. The affine combination $\alpha R+\beta S$ has barycentric coordinates $\alpha x+\beta y$. Hence:
	$$\delta(\alpha R+\beta S)=(\alpha x+\beta y)^t\Delta (\alpha x+\beta y)$$
	
	Quadraticity condition (\ref{quadracon}) is then written as:
	$$\begin{aligned} (\alpha x+\beta y)^t&\Delta (\alpha x+\beta y)^t=\\&=\alpha(\alpha-\beta)x^t\Delta x+\beta(\beta-\alpha)y^t\Delta y+4\alpha\beta\left(\frac{x+y}{2}\right)^t\Delta\left(\frac{x+y}{2}\right)\end{aligned}$$
	which is now straightforward using linearity of matrix product.
	%
	%
\end{proof}
In other words, this theorem proves that the definition \ref{definequadratico} identifies precisely the family of  mappings that can be written in coordinates with the classical quadratic matrix product.

For the particular case $\affi{B}=\RR$ we conclude that the space $\Quad{(\affi{A},\RR)}=\Quad{\affi{A}}$ of quadratic functions is a vector space linearly identified with the space of symmetric $(n+1)$-square matrices: $\Quad{\affi{A}}\simeq \Sym_{n+1}(\RR)$ (using a referential). Hence, it is a vector space with dimension $\binom{n+2}{2}$.

\end{document}